\documentclass[a4paper,11pt,reqno]{amsart}

\usepackage{amssymb, comment}
\usepackage{mathrsfs,wasysym}
\usepackage{mathscinet}
\usepackage{hyperref}
\usepackage{graphicx}
\usepackage{verbatim}
\usepackage[margin=1in]{geometry}
\DeclareGraphicsExtensions{.pdf}
\usepackage{setspace}
\usepackage{color}
\definecolor{gray}{gray}{0.4}

\usepackage{tikz}
\usetikzlibrary{arrows, automata,positioning,calc,shapes,decorations.pathreplacing,decorations.markings,shapes.misc,petri,topaths}
\usepackage{tkz-berge}
  \usepackage{pgfplots}
  \pgfplotsset{compat=newest}
  \usetikzlibrary{plotmarks}
  \usepackage{grffile}
\newlength\figureheight
  \newlength\figurewidth
\pgfplotsset{%
    tick label style={font=\scriptsize},
    label style={font=\footnotesize},
    legend style={font=\footnotesize},
   	every axis plot/.append style={very thick}
}

\hypersetup{
    pdftitle={Changepoint}}
\hypersetup{pdfsubject={Privacy, Security, Probability Theory}}
\hypersetup{pdfauthor={Michel Mandjes}}

\newcommand{\pp}{\mathbb P}

\newcommand{\bs}{\boldsymbol}

\newcommand{\bsx}{{\bs x}}
\newcommand{\bsX}{{\bs X}}
\newcommand{\bsy}{{\bs y}}
\newcommand{\bsm}{{\bs \mu}}
\newcommand{\bsn}{{\bs \nu}}
\newcommand{\TT}{^{{\rm T}}}

\newtheorem{lem}{Lemma}

\theoremstyle{remark}

\newcommand{\vr}{\varrho}

\newcommand{\ve}{\varepsilon}
\newcommand{\vb}{\vspace{3.2mm}}

\begin{document}

\title
{Changepoint detection for dependent Gaussian sequences}


\date{\today}

\begin{abstract}In this paper easily applicable techniques are devised for detecting changepoints in autocorrelated Gaussian sequences. Our method proceeds by sequential evaluation of a {\sc cusum}-type test statistic, which is compared to a predefined threshold. We assume that data is tested in sliding windows of fixed size. The distinguishing feature of this work is that, based on large deviations theory, we derive rather explicit equations that determine the threshold in such a way that the false alarm probability per window is approximately kept at the desired level. This criterion -- as opposed to the usual average run length -- allows to restrict not only the average number of false alarms but also their variability. Illustrative examples are provided, including the detection of a shift in mean in {\sc arma} processes. The procedures are validated by means of a broad set of simulation experiments, and overall perform well.
\end{abstract}

\author{W.\ Ellens$\,^{\bullet,\star}$, J.\ Kuhn$\,^{\bullet,\star}$, M.\ Mandjes$\,^{\bullet}$, P. $\dot{\rm Z}$uraniewski$\,^{\star,\circ}$}

\date{\today}

\maketitle

\noindent $^\bullet$ Korteweg-de Vries Institute for Mathematics,
University of Amsterdam, the Netherlands. 

\noindent $^\star$ TNO Performance of Networks and Systems, Delft, the Netherlands.

\noindent $^\circ$ Department of Applied
Mathematics, AGH University of Science and Technology, Krak\'ow,
Poland

\noindent M.\ Mandjes is also with  E{\sc urandom}
(Eindhoven University of Technology, the Netherlands) and  CWI
(Amsterdam, the Netherlands).

\noindent Email: {\scriptsize \url{wendy.ellens@tno.nl, j.kuhn@uva.nl, m.r.h.mandjes@uva.nl, piotr.zuraniewski@tno.nl}}

\section{Introduction}
The ability to detect \emph{changepoints} in data sequences (corresponding to a change in the underlying probability distribution) is of great practical importance, and one of the main concerns in statistical process control (SPC). In numerous application domains one is faced with problems of this nature. To mention but a few examples, changepoint
techniques have been used in finance \cite{CHEN}, electrocardiogram analysis \cite{GUST, GUST2}, climate change \cite{BARN} and communication networks \cite{CALL, THOT}.  

In the basic changepoint detection problem the goal is to detect a changepoint in a sequence of independent observations of some quality variable of interest. For continuous data it is usually assumed that the data is independent and normally distributed \cite{STOU}, and the change of interest is often a shift in the mean value. The goal is to detect it as soon as possible, while at the same time limiting the number of false alarms. Theoretical background on changepoint detection can be found in the books \cite{BANI,POOR}. The survey \cite{STOU} gives an overview on SPC from a practical perspective.

A commonly used technique in changepoint detection is that of {\it Cumulative Sum} ({\sc cusum}) \cite{PAGE}. For independent data, the {\sc cusum} statistic for detecting a change from the in-control parameter value to a pre-specified alternative can be expressed in terms of cumulative sums of log-likelihood ratio ({\sc llr}) increments. The monitoring is stopped and an alarm is issued as soon as the test statistic exceeds some predefined threshold. 

In the literature the question of how the threshold should be chosen is often disregarded, and when it is not, then typically the threshold has been selected based on a condition on the {\it average run length} ({\sc arl}), the expected time till the first false alarm \cite{BANI}. Since the {\sc arl} is simply the average of the stopping time, an obvious drawback of this approach is that it does not allow to restrict the variability of the false alarms. This can be crucial when thinking of applications in networks: Imagine, for example, one were to monitor patients' health data in a hospital (thus, testing multiple independent data streams in parallel). Then a high variance of false alarms could lead to a scenario where the capacity of the hospital staff is exceeded because a large number of false alarms (next to actual ones) occurred at the same time. 

Therefore, instead of the {\sc arl}, in this paper we restrict the probability of raising a false alarm in any given window. This ensures that the false alarm probability is low locally  (and is thus also still low on average). Furthermore, this approach circumvents an issue pointed out in \cite{Mei2008}, namely that the {\sc arl} is not always finite, and thus not in general an informative criterion. 

In his influential paper \cite{Lai1998}, Lai proposed two other false alarm criteria as an alternative to the {\sc arl}, which allow to limit the variance of the number of false alarms. It turns out that asymptotically, for window-limited detection, our criterion and one of Lai's are similar (see Section \ref{S4}). With respect to Lai's second criterion (which was coined {\it maximum local false alarm probability} in \cite{Tartakovsky2014book}), our method has the advantage of simplicity. In fact, it was stated in \cite[Ch.\ 8]{Tartakovsky2014book} that the practical implementation of this criterion is difficult because no closed-form expressions or even bounds are available that would allow the selection of the threshold. 

Previous results on how to select the threshold usually restrict the data points to be independent. For example, under this assumption the conceivable fact is proven that (under an appropriate scaling) 
a functional central limit theorem ({\sc clt}) holds, meaning that the cumulative random walk process converges to a Brownian motion. This result enables us to assess the test's false alarm probability \cite{SIEG}.  Apart from the {\sc clt} regime, asymptotic expressions for the false alarm probability have been derived under a large deviations scaling as well, see e.g.\  \cite[Ch. VI.E]{BUCK} and \cite{DESH}. 
Because these asymptotic expressions are available in closed form, choosing the threshold based on these results is relatively easy, yet ensures that the false alarm probability is limited. 


\vb

The analysis complicates significantly, however, if the observations do {\it not} correspond to independent variables. This situation is highly relevant, as in many practical situations the observations constituting the data sequence cannot be assumed independent. In the networking context, we refer to, e.g., the nice (unpublished) overview \cite{WIL} for an extensive treatment of traffic characteristics in communication networks; notably, it has been found that there are non-negligible correlations over broad ranges of time scales.

This motivates that in the current paper we focus on Gaussian processes that exhibit serial dependence. 
An important class of Gaussian processes that include dependence is that of the so-called 
autoregressive moving-average (short: {\sc arma}) processes \cite{BJ,BROC}, which we consider as a more specific example. 
For the class of {\sc arma} processes
Johnson and Bagshaw \cite{JOHN} established the convergence to Brownian motion, thus enabling the type I error (false alarm) analysis of a {\sc cusum}-type procedure. Alternative tests under the {\sc clt} scaling were described extensively by
Cz\"org\H{o} and Horv\'ath \cite[Ch. IV]{CZOR}, with a focus on a Brownian-bridge based test statistic
(see also \cite{HUSK}). Basseville and Nikiforov \cite[Ch.\ 7]{BANI} discuss testing procedures for dependent Gaussian processes that rely on a whitening transformation of the data sequence. A similar avenue is taken in \cite{ENGL} and \cite{ROBB} for the problem of mean shift detection in {\sc arma} processes. Besides these works, upper bounds have been provided for more general scenarios, where the Gaussianity assumption is relaxed (see e.g.\ \cite{Lai1998} and \cite{TartakovskyRBK2006}).

\vb

The current paper contributes to the theory on changepoint techniques for serially correlated data. We develop a window-limited testing procedure with {\sc llr} test statistic (in the spirit of the {\sc cusum} method), and provide a method for selecting the threshold (function) such that the probability of raising a false alarm is low in every given window of data points, as motivated above. An advantage of testing data in windows rather than keeping the entire history of observations is that a change can be detected more quickly since it has a bigger impact relative to the (fewer) previous observations within the current window. Furthermore, the usual assumption of stationary data is less restrictive in this case.

While previous (asymptotic) work on {\sc cusum} for dependent data has primarily focused on the {\sc clt} regime, in the present paper we consider a large-deviations (short: {\sc ld}) setting. 
More specifically, we construct {\sc ld}-based {\sc cusum}-type changepoint detection tests for dependent normal data, covering also the class of (Gaussian) {\sc arma} processes. 
Since {\sc ld}  theory \cite{BUCK,DEMB}
focuses on the rare-event setting, this framework is particularly suitable for the problem at hand as the probability of raising a false alarm is required to be low.  

An additional attractive feature of applying {\sc ld} theory here is that it nicely facilitates the analysis of hypothesis testing with multiple alternatives. In the changepoint detection problem we have to consider a union of hypotheses corresponding to a change in a parameter value {\it at some point} in the dataset. In the {\sc ld} regime the probability of such a union of events essentially coincides with the probability of the most likely event among them; this phenomenon is usually referred to as the {\it principle of the largest term} \cite{GOW}.  We therefore obtain a threshold {\it function} rather than a single value as is usually assumed (see \cite{BANI,Tartakovsky2014book}), ensuring that the probability of raising a false alarm is essentially equally likely irrespective of the location of the changepoint. We provide a numerical example in Section \ref{S5} that indicates that choosing a threshold \emph{function} is indeed favourable.

In that section, we also discuss a number of relevant cases in greater detail: a change in the mean (with the correlations held fixed), a change in variance (for independent observations), and a change of the `scale' of the process (that is, the means blow up by a factor $f$, the covariance matrix by a factor $f^2$). In these cases we obtain particularly simple equations for the threshold function, see Eqs. (\ref{beetje}), (\ref{b_variance}) and (\ref{b_scale}), respectively. The change in scale example was considered in more detail in \cite{KEM2014} in a multidimensional setting; it has applications in the context of communication networks where a change in scale may result from an increase of the number of users.

The paper is organized as follows. In Section \ref{S2} we provide preliminaries on {\sc cusum}, reviewing the {\it independent case} in the {\sc ld} scaling. Then Section \ref{S3} provides a series of useful computations for likelihood ratio tests related to multivariate normal distributions, which are used in Section \ref{S4} to develop changepoint detection tests for {\it dependent} data, and includes the aforementioned more specific examples. Section \ref{S5} presents an extensive simulation study so as to assess the performance of the tests; these experiments confirm that the proposed procedure works well in a broad range scenarios.

\section{Cumulative Sum: preliminaries}\label{S2}
Consider a representative window of observations $X_1,X_2,\ldots, X_n$, during which potentially a changepoint occurs. In this section we assume that the $X_i$ are independent, but we do not assume anything about their distribution. Later in this paper we look at situations in which the $X_i$ may be dependent, but follow a normal distribution. In probabilistic terms a {\it changepoint}, to be considered as a change in the statistical law of the
underlying random variable, can be described as follows.
\begin{itemize}
\item Under the null-hypothesis (H$_0$) the $X_i$ ($i=1,\ldots,n$)
are independent and identically distributed (i.i.d.) realizations of a random variable with density $f(\cdot)$.
\item Under the alternative hypothesis (H$_1$) up to $k-1$ the observations are i.i.d.\ samples from a distribution
with density $f(\cdot)$, while from observation $k$ on they are i.i.d.\ with a {\it different} density $g(\cdot)$ (for some $k$ ranging between $1$ and $n$).
\end{itemize}
In other words: under the null-hypothesis there has {\it not} been a changepoint, while under the alternative hypothesis the process changes. Observe that this setup is not a simple binary hypothesis testing problem, as the alternative is essentially a {\it union} of hypotheses. More precisely: with H$_1(k)$ corresponds to having a changepoint at $k$, we can write H$_1$
as the union of the H$_1(k)$, with $k=1,\ldots,n.$

A changepoint detection test, that is, a test that determines whether to accept the null hypothesis or to reject it ---~in which case it issues an alarm~---- aims at keeping the probability of a type I error (a false alarm) limited. On the other hand, the test should be such that the detection probability is as high as possible, in other words, it should minimize the type II error probability while maintaining the false alarm rate at a given low level. 

The technique we describe in this section, known as {\sc cusum}, has been proposed \cite{PAGE} to identify parameter changes from the in-control value to a pre-specified alternative. Since in practice the parameter after the change is typically unknown, it is commonly replaced by its maximum likelihood estimator (resulting in the generalized likelihood ratio ({\sc glr}) test), or by some smallest tolerable value \cite{BASS}. Also a combination of multiple testing procedures is possible, as, for example, proposed in \cite{ZHAO}. Since this question is not in the scope of the current paper, in the following description of the {\sc cusum} method we assume that the alternative is specified ---~we roughly follow the setup presented in \cite[Ch.\ II.6]{SIEG}.

Consider first the common likelihood test for H$_0$ versus H$_1(k)$. Evidently, the statistic to be considered is
\[\bar S_k:=\left(\prod_{i=k}^n g(X_i)\right)\left/\left(\prod_{i=k}^n f(X_i)\right)\right.\hspace{-1mm};\]
it turns out, though, that it is more practical to work with the corresponding {\it log}-likelihood:
\[S_k :=\sum_{i=k}^n \log\left(\frac{g(X_i)}{f(X_i)}\right).\]
To deal with the fact that H$_1$
equals the union of the H$_1(k)$,
we have to verify whether there is a $k\in\{1,\ldots,n\}$ such that $S_k$ exceeds a certain critical value.
As a result, the statistic for the {\it composite} test (that is,  H$_0$ versus H$_1$) is
\begin{equation}
\label{rp}t_n:=\max_{k\in\{1,\ldots,n\}} S_k=T_n-\min_{k\in\{1,\ldots,n\}}T_{k-1},\end{equation}
with $T_k$ denoting the cumulative sum $\sum_{i=1}^{k} \log\big[ g(X_i)/f(X_i)\big];$
the null-hypothesis is rejected if $t_n$ exceeds some critical level $b$.

Observe from the above that the test statistic can
be written in terms of the cumulative sums $T_k$ (corresponding to increments that are distributed
as $g(X_i)/f(X_i)$), which explains the name of the test. Also, note that the statistic (\ref{rp}) represents
the height of the random walk $T_k$ relative to the minimum that was achieved so far; in this sense, there
is a close connection to an associated (discrete-time)  {\it queueing} process, as described in, e.g., \cite{SIEG}. {\sc cusum} has certain optimality problems in terms of the tradeoff mentioned above (timely detection versus low rate of false alarms, that is), as established  in a Bayesian framework in \cite{SHIR1, SHIR2}, whereas \cite{LORD,POLL} address this property in the non-Bayesian setting.

We now scale the threshold $b$ by $n$, and focus on asymptotics for large $n$; this limiting regime is usually referred to as the {\it large deviations regime} \cite{BUCK, DEMB,MAND}. More specifically, we analyze the probability of issuing a false alarm (type I error), that is, ${\mathbb P}_0(t_n\ge nb)$. Here ${\mathbb P}_0$ corresponds to probability under H$_0$ and ${\mathbb E}_0$ is the associated expectation. We roughly follow the setup of \cite[Ch. VI.E]{BUCK}.
Under H$_0$, due to reversibility arguments,
\[t_n=T_n-\min_{k\in\{1,\ldots,n\}}T_{k-1} =
\max_{k\in\{1,\ldots,n\}}(T_n-T_{k-1})\stackrel{\rm d}{=}  \max_{k\in\{1,\ldots,n\}}T_k,\]
so that the probability of our interest can be rewritten as
\[{\mathbb P}_0(t_n\ge nb) = {\mathbb P}_0\left(\exists k \in \{1,\ldots,n\}: T_k \ge nb\right).\]
Due to $n^{-1}\cdot \log n\to 0$ and
\[\max_{k \in \{1,\ldots,n\}}{\mathbb P}_0\left(T_k \ge nb\right)\le
{\mathbb P}_0\left(\exists k \in \{1,\ldots,n\}: T_k \ge nb\right)\le
n\cdot \max_{k \in \{1,\ldots,n\}}{\mathbb P}_0\left(T_k \ge nb\right),\]
we have the following expression for the so-called \emph{decay rate}
\[\lim_{n\to\infty}\frac{1}{n}\log {\mathbb P}_0(t_n\ge nb) = \max_{\lambda\in(0,1]}
\lim_{n\to\infty}\frac{1}{n}\log {\mathbb P}_0\left(\frac{T_{n\lambda}}{n}\ge b\right)\]
(realize that $n\lambda$ is not necessarily integer, so there is mild abuse of notation in the previous display); in words, this means that the decay rate of the union of all $n$ events coincides with the decay rate of the most likely event among these (the so-called  `principle of the largest term'; see \cite{GOW}). 
Relying on Cram\'er's theorem \cite[Ch. II.A]{BUCK}, we can rewrite the above decay rate to
\[\lim_{n\to\infty}\frac{1}{n}\log {\mathbb P}_0(t_n\ge nb) = \max_{\lambda\in(0,1]}
\lim_{n\to\infty}\frac{\lambda}{n\lambda}\log {\mathbb P}_0\left(\frac{T_{n\lambda}}{n\lambda}\ge \frac{b}{\lambda}\right)=
\max_{\lambda\in(0,1]}
\left(-\lambda\sup_{\theta}\left(\theta\frac{b}{\lambda}-\log M(\theta)\right)\right);\]
here $M(\theta)$ is the moment generating function (under H$_0$) of $\log\big[ g(X_i)/f(X_i)\big]:$
\[M(\theta)= {\mathbb E}_0 \exp\left({\theta \log\frac{ g(X_i)}{f(X_i)}}\right) = {\mathbb E}_0 \left(\frac{g(X_i)}{f(X_i)}\right)^\theta =
\int_{-\infty}^\infty (g(x))^\theta (f(x))^{1-\theta}{\rm d}x.\]
We can then set $b$ such that the decay rate under study equals some predefined (negative) constant $-\gamma$ (where $\gamma>0$).
In principle, however, there is no need to take a {\it constant} $b$; we could pick a {\it function} $b(\lambda)$ instead. It can be seen that, in terms of optimizing the type II error performance, it is optimal to choose this function $b(\lambda)$ such that
$$\lim_{n\to\infty}\frac{1}{n}\log {\mathbb P}_0\left(\frac{T_{n\lambda}}{n}\ge b(\lambda)\right)=-\lambda\sup_{\theta}\left(\theta\frac{b(\lambda)}{\lambda}-\log M(\theta)\right)$$
is {constant} in $\lambda\in(0,1]$ (and equaling $-\gamma$). Intuitively, this choice entails that for any point $n\lambda$ in time, issuing an alarm (which is done if $T_n-T_{n\lambda-1}$ exceeds $nb(1-\lambda+1/n)$) is essentially equally likely if there is no changepoint.

In the setup described above the individual observations $X_i$ are assumed to be {independent}. The main objective of the paper is to develop a machinery that can deal with {\it dependent} data. As mentioned earlier, we focus on the case that the data stem from a multivariate normal distribution. To this end, we first work out the likelihood ratio test of a single multivariate normal distribution against another one in Section \ref{S3}, which is used in Section \ref{S4} to develop a changepoint detection procedure for dependent normal data.

\section{Likelihood ratio test for multivariate normal data}\label{S3}

As we saw in the previous section, the {\sc cusum} method is in essence a sequentially applied {\sc llr} hypothesis test. We therefore first consider the situation that under H$_0$ the data $X_1,\ldots,X_n$ has a normal distribution with mean $\bar\mu$ under $H_0$ and mean $\bar\nu$ under $H_1$.
That is, in this section we assume that there is no changepoint (or, equivalently, that the change has occurred already at the first observation within the considered window). The results of this section will be used in Section \ref{S4} to develop a procedure to find a change \emph{somewhere} in the sequence.

It is immediately seen that, without loss of generality, we can pick $\bar\mu=0$ (by subtracting it from $\bar\nu, X_1,\ldots,X_n$). Because we wish to explicitly allow for correlated data points, we further assume that the vector of observations ${\bs X}=(X_1,\ldots,X_n)$ stems from an $n$-dimensional multivariate normal distribution with mean $\bs\mu_n\equiv\bs\mu$ and covariance matrix $\Sigma_n\equiv \Sigma$ (which thus does not need to be diagonal), denoted by ${\mathscr N}({\bs \mu},\Sigma)$, while under H$_1$ they stem from ${\mathscr N}({\bs \nu},T)$.
 
We let $f_n(\cdot)$ and $g_n(\cdot)$ be the corresponding $n$-dimensional densities, that is,
\[f_n(\bsx)=(2\pi)^{-n/2} \,|\,\Sigma\,|^{-1/2}\ \exp\left(-\frac{1}{2}\bsx \TT \Sigma^{-1} \bsx\right),\]
and
\[g_n(\bsx)=
(2\pi)^{-n/2} \,|\,T\,|^{-1/2}\exp\left(-\frac{1}{2}(\bsx-\bsn) \TT T^{-1} (\bsx-\bsn)\right).\]
Observe that $\bsm$ and $\bsn\in {\mathbb R}^n$, while
$\Sigma$ and $T$ are positive-definite matrices of dimension $n\times n$.
In this section, we first develop a large-deviations based likelihood ratio test for distinguishing $g_n(\cdot)$ from $f_n(\cdot)$, and then specialize to a series of relevant special cases.

A {\sc llr} hypothesis test features the test statistic
 \[{\mathscr L}_n(\bsX)= \log\left(\frac{ g_n({\bs X})}{f_n({\bs X})}\right),\] which
can be evaluated as
\begin{equation}\label{LL}{\mathscr L}_n(\bsX)=\frac{1}{2}\log|\,\Sigma\,|-\frac{1}{2}\log|\,T\,|+\frac{1}{2}\bsX\TT\Sigma^{-1}\bsX-\frac{1}{2}(\bsX-\bsn)\TT T^{-1}(\bsX-\bsn).\end{equation}


To determine the critical value $nb$ above which the null hypothesis is rejected, we wish to evaluate the type I error probability ${\mathbb P}_0({\mathscr L}_n(\bsX) \ge nb)$, where $b>{\mathbb E}_0 {\mathscr L}_n(\bsX) /n$.
It turns out to be hard to evaluate this probability explicitly, but we can derive an accurate approximation based on large deviations theory.
Relying on the G\"artner-Ellis theorem \cite{BUCK,DEMB}, the following equation holds for the decay rate
\[\lim_{n\to\infty}\frac{1}{n}\log
{\mathbb P}_0({\mathscr L}_n(\bsX) \ge nb) =
-{\mathscr I}(b),\]
where ${\mathscr I}(b)$ denotes the associated Legendre transform
\begin{equation}\label{LT}{\mathscr I}(b):=\sup_{\theta}\left(\theta b -\lim_{n\to\infty}\frac{1}{n}
\log {\mathbb E}_0 \exp (\theta {\mathscr L}_n(\bsX))\right),\end{equation}
given that the limiting log-moment generating function exists.
This leads to the approximation
\[{\mathbb P}_0({\mathscr L}_n(\bsX) \ge nb) \approx
e^{-n{\mathscr I}(b)}.\]
To use this approximation, we first compute the moment generating function ${\mathbb E}_0 \exp (\theta {\mathscr L}_n(\bsX))$
in more explicit terms.
It is clear that
\[{\mathbb E}_0 \exp (\theta {\mathscr L}_n(\bsX))=
(2\pi)^{-n/2} \,|\,\Sigma\,|^{-1/2}\int_{-\infty}^\infty
\cdots\int_{-\infty}^\infty  \exp (\theta {\mathscr L}_n(\bsx)) \exp\left(-\frac{1}{2}\bsx\TT\Sigma^{-1}\bsx\right){\rm d}x_1\cdots {\rm d}x_n.
 \]
 Then notice that
 \begin{equation}
 \label{fo}\theta {\mathscr L}_n(\bsx) -\frac{1}{2}\bsx\TT\Sigma^{-1}\bsx = \frac{\theta}{2}\log\frac{|\,\Sigma\,|}{|\,T\,|}-\frac{1}{2}
 \bsx\TT(\theta T^{-1}+(1-\theta)\Sigma^{-1})\bsx+\theta\bsn\TT T^{-1}{\bsx}-\frac{\theta}{2}\bsn\TT T^{-1}\bsn.\end{equation}
Now realize that $\theta T^{-1}+(1-\theta)\Sigma^{-1}$ is positive-definite; let $B\TT B$ be the corresponding Cholesky decomposition.
As a next step, we perform the substitution $\bsy = B\bsx$, so that
\[{\rm d}x_1\cdots {\rm d}x_n= |\,B^{-1}\,| \,{\rm d}y_1\cdots {\rm d}y_n=\frac{1}{|\,\theta T^{-1}+(1-\theta)\Sigma^{-1}\,|^{1/2}} \,{\rm d}y_1\cdots {\rm d}y_n.\]
Then Expression (\ref{fo}) can be rewritten as
\[\frac{\theta}{2}\log\frac{|\,\Sigma\,|}{|\,T\,|}-\frac{1}{2}\bsy\TT\bsy +\theta\bsn\TT T^{-1} B^{-1} \bsy-\frac{\theta}{2}\bsn\TT T^{-1}\bsn,\]
which equals
\begin{eqnarray*}
\lefteqn{\frac{\theta}{2}\log\frac{|\,\Sigma\,|}{|\,T\,|}-\frac{1}{2}\left(\bsy-\theta(B^{-1})\TT T^{-1}\bsn\right)\TT\left(\bsy-\theta(B^{-1})\TT T^{-1}\bsn\right)}\\
&-&\frac{\theta}{2}\bsn\TT T^{-1}\bsn
+\frac{\theta^2}{2}\bsn\TT T^{-1}\left(\theta T^{-1}+(1-\theta)\Sigma^{-1}\right)^{-1}T^{-1}\bsn.\end{eqnarray*}
Recognizing a multivariate normal density, we conclude that the moment generating function ${\mathbb E}_0 \exp (\theta {\mathscr L}_n(\bsX))$ equals, with $I_n$ denoting an $n\times n$ identity matrix,
\begin{eqnarray}
{\mathbb E}_0 \exp (\theta {\mathscr L}_n(\bsX))&=&{\left(\frac{|\,\Sigma\,|}{|\,T\,|}\right)^{\theta/2}\frac{|\,\Sigma\,|^{-1/2}}{|\,\theta T^{-1}+(1-\theta)\Sigma^{-1}\,|^{1/2}}}\nonumber\\
&&\times
\exp\left(-\frac{\theta}{2}\bsn\TT T^{-1}\bsn
+\frac{\theta^2}{2}\bsn\TT T^{-1}\left(\theta T^{-1}+(1-\theta)\Sigma^{-1}\right)^{-1}T^{-1}\bsn\right)\nonumber\\
&=&\left(\frac{|\,\Sigma\,|}{|\,T\,|}\right)^{\theta/2}\frac{1}{|\,\theta T^{-1}\Sigma+(1-\theta)I_n\,|^{1/2}}\label{momentgeneratingfunction}\\
&&\times
\exp\left(-\frac{\theta}{2}\bsn\TT T^{-1}\bsn
+\frac{\theta^2}{2}\bsn\TT T^{-1}\left(\theta T^{-1}+(1-\theta)\Sigma^{-1}\right)^{-1}T^{-1}\bsn\right).\nonumber\end{eqnarray}

The above analysis gives, in principle, a technique to calculate ${\mathscr I}(b)$, and hence, a technique
to approximate the type I error probability. This allows us to determine the critical value $b$. In specific cases, the computations can be made more explicit. Below
we treat two of those special cases. In Section \ref{sec_diff_mean} we work out the moment generating function (\ref{momentgeneratingfunction}) and find the Legendre transform (\ref{LT}) for a test designed to decide between two different means, while for the special case of independent data (\ref{momentgeneratingfunction}) is simplified in Section \ref{sec_diff_var} .

\subsection{Special case I: difference in mean for dependent data} \label{sec_diff_mean}
In the first special case we focus on, there is only a difference in the means of the multivariate normal distributions, that is, the covariance matrix is left unchanged: $\Sigma =T$. It means that
\[{\mathbb E}_0 \exp (\theta {\mathscr L}_n(\bsX))=\exp\left(-\frac{\theta}{2}\bsn\TT T^{-1}\bsn
+\frac{\theta^2}{2}\bsn\TT T^{-1}\bsn\right). \]
As a consequence ---~defining ${\mathscr I}_n(b):=n{\mathscr I}(b)$~--- we have
\[{\mathscr I}_n(b)= \sup_\theta\left(n\theta b +\frac{\theta}{2}\bsn\TT T^{-1}\bsn
-\frac{\theta^2}{2}\bsn\TT T^{-1}\bsn\right).\]
The supremum can be determined explicitly, leading to
\begin{equation}
\label{nun}{\mathscr I}_n(b)= \frac{(nb+\frac{1}{2}\bsn\TT T^{-1}\bsn)^2}{2\bsn\TT T^{-1}\bsn}.\end{equation}
We will use this result in Section \ref{SUB41} to develop a changepoint detection test to find a change in the mean of a dependent (multivariate normal) sequence.

\subsection{Special case II: difference in mean and variance for independent data}\label{sec_diff_var}
In the second special case we have that there is a difference in both mean and covariance matrix of
the multivariate normal distributions, but in such a way that the covariance matrices $\Sigma$ and $T$ correspond
to independent random variables.  In this setting $\Sigma$ is the diagonal matrix with the vector ${\bs \sigma}^2$ on the diagonal (to be denoted by
$ {\rm diag}({{\bs \sigma}}^2)$), while $T={\rm diag}({\bs \tau}^2).$
It is a matter of elementary calculus to verify that
\begin{eqnarray}
{{\mathbb E}_0 \exp (\theta {\mathscr L}_n(\bsX))}&=&{
\prod_{i=1}^n \left(\frac{\sigma_i}{\tau_i}\right)^\theta\times
\prod_{i=1}^n\left(\theta\frac{\sigma_i^2}{\tau_i^2}+(1-\theta)\right)^{-1/2}}\nonumber\\
&&\times \exp\left(-\frac{\theta}{2}\sum_{i=1}^n \frac{\nu_i^2}{\tau_i^2}+\frac{\theta^2}{2}
\sum_{i=1}^n \frac{\nu_i^2 \sigma_i^2/\tau_i^2}{\theta \sigma_i^2+(1-\theta)\tau_i^2}\right).\label{mom_gen_fun_indep}
\end{eqnarray}
\noindent
The above result is used in Section \ref{change_variance_independent} for a test that detects a change in variance somewhere in a sequence of independent normally distributed data.

\section{Changepoint detection tests for dependent data}\label{S4}
We now propose a series of changepoint detection tests, in line with the one presented for an i.i.d. sequence in \cite[Ch.\ VI.E]{BUCK} (discussed in Section \ref{S2} of this paper).
The idea is that H$_0$ corresponds to a model ${\mathbb P}_0$, whereas under H$_1$ there is a shift of the model
${\mathbb P}_0$ to ${\mathbb P}_1$ at the $(n\beta+1)$-th observation, for some $\beta\in [0,1)$ such that $n\beta$ is integer-valued. In line with \cite[Ch.\ VI.E, Eq. (43)]{BUCK}
we reject H$_0$ if
\begin{equation}
\label{testcrit} \max_{\beta\in[0,1)} \left( \frac{1}{n} {\mathscr L}_{n,\beta}(\bsX) -{b}(\beta) \right) := \max_{\beta\in[0,1)} \left( \frac{1}{n}\log \left(\frac{g_{n,\beta}({\bs X})}{f_n({\bs X})}\right) -{b}(\beta) \right) > 0,
\end{equation}
where the density $g_{n,\beta}(\cdot)$ corresponds to H$_1$ with a change at time $n\beta+1$, and ${b}(\cdot)$ is a function specified below.
Large-deviations theory enables us to compute
\[\lim_{n\to\infty}\frac{1}{n}\log {\mathbb P}_0\left(\max_{\beta\in[0,1)} \big( {\mathscr L}_{n,\beta}(\bsX) -{b}(\beta) \big) > 0\right),\]
using the machinery of Section \ref{S3}. To optimize the type II error rate performance \cite[Ch.\ VI.E, p.\ 113]{BUCK}, ${b}(\cdot)$ should be chosen such that the decay rate satisfies
\begin{equation}
\label{unif}-{\mathscr I}(b(\beta))=\lim_{n\to\infty}\frac{1}{n}\log {\mathbb P}_0\big(  {\mathscr L}_{n,\beta}(\bsX) -{b}(\beta) > 0\big)=-\gamma\end{equation}
for a uniform positive $\gamma$, across all $\beta\in[0,1)$; this  enables us to determine ${b}(\beta).$ In practice the observations arrive one by one and at every new observation $X_m$ the changepoint detection test is then performed on the sequence of the $n$ most recent observations $(X_{m-n+1},\ldots X_m)$. An alarm is issued at time $m$ if the test statistic ${\mathscr L}_{n,\beta}(\bsX)$ exceeds the threshold ${b}(\beta)$ for any $\beta\in[0,1)$. The goal is to detect a changepoint as soon as possible, while at the same time keeping the number of false alarms limited. We explain the details of how to carry out the testing in more detail in the remainder of this section, and also provide numerical examples in Section 5. In the following, we again use $X_1,\dots,X_n$ to represent the observations of the current window (thus, dropping the enumeration of the windows by $m$).

\vb 

Note that traditionally in changepoint detection the {\sc arl} --- the expected time till the first false alarm --- has been used to design procedures with a limited number of false alarms. However,  the {\sc arl} criterion does not ensure that the number of false alarms is low for every window (see also the discussion in \cite{Lai1998}), and furthermore, it may not always be applicable (see \cite{Mei2008} for examples where the {\sc arl} becomes infinite). Our methodology in the current paper circumvents these issues.

It should be pointed out that the proposed procedure is essentially a (window-limited) CUSUM-type detection procedure. In \cite{Lai1998} CUSUM is proven to be \emph{asymptotically optimal} (as $\alpha\to 0$) in that it achieves the lowest possible detection delay provided that $\sup_{m\geq 1} \pp_0(m\leq\tau\leq m+n)\leq \alpha$, where $\tau$ denotes the stopping time of the CUSUM-type procedure. It turns out that the distribution of $\tau$ is approximately exponential \cite[Ch.\ 8]{Tartakovsky2014book}. Therefore, from the memoryless-property of the exponential distribution, we have that 
\begin{align*}
\sup_{m\geq 1} \pp_0(m\leq\tau\leq m-1+n)&\approx \pp_0(1\leq \tau\leq n)\,.
\end{align*}
Since we consider a window-limited detection procedure, where $\tau<n$ is not considered, we thus impose (\ref{unif}) rather than the criterion proposed by Lai in \cite{Lai1998}.
Furthermore, as (\ref{unif}) limits the false alarm probability for any given window, the resulting \emph{average} false alarm rate (averaged over all windows that do not include the changepoint) will also be limited to the same level.

\vb

We now perform the computation of (\ref{unif}) and the determination of the critical function $b(\beta)$ for various specific models. In \cite[Ch.\ VI.E Example 3]{BUCK} the critical function is determined for a change in mean in a sequence of independent normally distributed observations. In Section \ref{SUB41} we look at a change in mean somewhere in a (dependent) multivariate normal sequence (using the result of Section \ref{sec_diff_mean}), in Section \ref{change_variance_independent} we consider a change in variance for independent normally distributed sequences (using the result of Section \ref{sec_diff_var}) and Section \ref{change_scale} treats the case of a change in scale of a (dependent) multivariate normal sequence.

\subsection{Test 1: change in mean for dependent data}\label{SUB41}

In this section we show how to compute the critical function $b(\beta)$ when testing for a change in the mean of a dependent sequence. We derive an explicit expression for $b(\beta)$ for the case of autoregressive-moving-average ({\sc arma}) processes.

We are in the setting that $\Sigma=T$, and that we want to detect a change in mean at some
index $n\beta+1$, for $\beta\in[0,1).$ Without loss of generality we consider a change from mean 0 to some other value, say $\bar\nu$. In line with the above, we wish to find a function ${b}(\beta)$ such that
(\ref{unif}) holds for $\beta\in[0,1)$,
for a given $\gamma>0.$ We can apply formula (\ref{nun}), with the first $n\beta$ entries of $\bsn$ equal to 0 and the last $n(1-\beta)$ equal to $\bar\nu.$ Defining
\[t_{n,\beta}:=\sum_{i=n\beta+1}^n\sum_{j=n\beta+1}^n (T^{-1})_{i,j},\]
we obtain
\[-\gamma=\lim_{n\to\infty}\frac{1}{n}\log {\mathbb P}_0\left(\frac{1}{n}\log\left(\frac{g_{n,\beta}(\bsX)}{f(\bsX)}\right) \ge {b}(\beta)\right)=-{\mathscr I}(b(\beta))=
-\lim_{n\to\infty}\frac{1}{2}\frac{(n{b}(\beta)+\frac{1}{2}\bar\nu^2t_{n,\beta})^2}{n\bar\nu^2t_{n,\beta}}.\]
As an example we could consider $\bsX$ corresponding to an {\it autoregressive process of order 1} (usually abbreviated to {\sc ar}(1)). This is a stationary process (with mean $c$) obeying the recursion
\[X_i - c = \vr (X_{i-1}-c) + \ve_i,\]
where the $\ve_i$\,s are i.i.d.\ samples from a zero-mean normal distribution with variance $\sigma^2$ (where we assume $|\varrho|<1$). It is known that
\[T=\frac{\sigma^2}{1-\vr^2}
\left(\begin{array}{cccccc}
1&\vr&\vr^2&\vr^3&\cdots&\vr^{n-1}\\
\vr&1&\vr&\vr^2&\cdots&\vr^{n-2}\\
\vr^2&\vr&1&\vr&\cdots&\vr^{n-3}\\
\vr^3&\vr^2&\vr&1&\cdots&\vr^{n-4}\\
\vdots&\vdots&\vdots&\vdots&\ddots&\vdots\\
\vr^{n-1}&\vr^{n-2}&\vr^{n-3}&\vr^{n-4}&\cdots&1
\end{array}\right).\]
It is elementary to verify that
\[T^{-1}=\frac{1}{\sigma^2}\left(\begin{array}{cccccc}
1&-\vr&0&0&\cdots&0\\
-\vr&1+\vr^2&-\vr&0&\cdots&0\\
0&-\vr&1+\vr^2&-\vr&\cdots&0\\
0&0&-\vr&1+\vr^2&\cdots&0\\
\vdots&\vdots&\vdots&\vdots&\ddots&\vdots\\
0&0&0&0&\cdots&1
\end{array}\right).\]
It follows that (realizing that there are roughly $n$ diagonal entries of value $1+\vr^2$, and that there are roughly
$2n$ entries of value $-\vr$ above and below the diagonal),
\[\lim_{n\to\infty}\frac{t_{n,\beta}}{n(1-\beta)}=\frac{1}{\sigma^2}\left(1\cdot ({1+\vr^2})+2\cdot (-\vr)\right)=\left(\frac{1-\vr}{\sigma}\right)^2,\]
and hence
\begin{equation}\label{crit_func_AR1}
{b}(\beta)=\bar\nu\left(\frac{1-\vr}{\sigma}\right)\sqrt{2\gamma(1-\beta)}-\frac{1}{2}\bar\nu^2\left(\frac{1-\vr}{\sigma}\right)^2(1-\beta).\end{equation}
Compared to the function ${b}(\beta)$ that was derived
for the unit-variance i.i.d.\ case \cite[Ch.\ VI.E, p.\ 113]{BUCK}, $\bar{\nu}$ needs to be replaced by $\bar\nu(1-\vr)/\sigma$, in order to account for the dependence between the observations, and the value of the variance. For $\vr=0$ and $\sigma^2=1$, the two functions obviously match.

\vb

Also in case that $T^{-1}$ cannot be computed explicitly, we can still find the limiting value of $t_{n,\beta}/(n(1-\beta)).$ We now consider the general {\sc arma}($p,q$) model, defined as a stationary model with mean value $c$ obeying \cite{BROC}
\begin{equation}\label{RE}X_i -c = \ve_i +\sum_{j=1}^p \vr_{j}(X_{i-j}-c) +\sum_{j=1}^q \vartheta_j\ve_{i-j},\end{equation}
for $p,q\in {\mathbb N}$, where we assume that the roots of the {\sc ar} polynomial lie outside the unit circle. Again we assume that the $\ve_i$ are i.i.d.\ samples from a zero-mean normal distribution with variance $\sigma^2.$

The following lemma implies that the limiting value of $t_{n,\beta}/(n(1-\beta))$
does not depend on $\beta$, or, put differently, that $t_{n,\beta}$ grows essentially linear in $n(1-\beta)$; cf.\ \cite[Eq.\ (9)]{ROBB}. 

\begin{lem}\label{lemma1} For ${\bs X}$ obeying an {\sc arma}($p,q$) model, and $\beta \in[0,1) $,
\[{\mathscr T}_\beta:=\lim_{n\to\infty}\frac{t_{n,\beta}}{n(1-\beta)}=\left(\frac{1-\sum_{j=1}^p\vr_j}{\sigma\left(1+\sum_{j=1}^q \vartheta_j\right)}\right)^2=:{\mathscr T}.\]
\end{lem}

The proof can be found in Appendix \ref{proof_lemma1}. The immediate consequence of the lemma is that
\[-\gamma=\lim_{n\to\infty}\frac{1}{n}\log {\mathbb P}_0\left(\frac{1}{n}\log\left(\frac{g_{n,\beta}(\bsX)}{f(\bsX)}\right) \ge {b}(\beta)\right)=-\frac{1}{2}\frac{\left(b(\beta)+\frac{1}{2}\bar{\nu}^2{\mathscr T}(1-\beta)\right)^2}{\bar\nu^2{\mathscr T}(1-\beta)},\]
and
\begin{equation}
\label{beetje}b(\beta)=\bar \nu\sqrt{2{\mathscr T}\gamma(1-\beta)}-\frac{1}{2}\bar\nu^2 {\mathscr T}(1-\beta).\end{equation}

We have seen that for {\sc ar}(1) processes ${\mathscr T}=((1-\vr)/\sigma)^2$. From Lemma \ref{lemma1} it follows that for an {\sc ma}(1) process with parameter $\vartheta$ it holds that ${\mathscr T}=1/(\sigma(1+\vartheta))^2$ and
\begin{equation}\label{crit_func_MA1}
{b}(\beta)=\bar\nu\left(\frac{1}{\sigma(1+\vartheta)}\right)\sqrt{2\gamma(1-\beta)}-\frac{1}{2}\bar\nu^2\left(\frac{1}{\sigma(1+\vartheta)}\right)^2(1-\beta).\end{equation}

\subsection{Test 2: change in variance for independent data}\label{change_variance_independent}
We now consider the case in which there is no change in mean, where under H$_0$ all observations are independent and normally distributed with variance $\sigma^2$ while under H$_1$ the variance changes from $\sigma^2$ to $\tau^2$ at some specific moment. We set $\bsn={\bs 0}$,
$\Sigma = \sigma^2I_n$, and $T$ is an $n\times n$ diagonal matrix with
$\sigma^2$ at the first $m=\beta n$ diagonal positions ($\beta\in[0,1)$), and $\tau^2$
at the other diagonal positions. Note that this corresponds to a change in variance at time $\beta n+1$. Filling out (\ref{mom_gen_fun_indep}), we get
\begin{eqnarray*}
\Lambda_\beta(\theta) &:=& \frac{1}{n}\log {\mathbb E}_0 \exp (\theta {\mathscr L}_n(\bsX))\\
&=& \theta (1-\beta)\log\frac{\sigma}{\tau}+\frac{1}{2}(1-\beta)\log\tau^2-\frac{1}{2} (1-\beta)\log\left(\theta{\sigma^2}+(1-\theta){\tau^2}\right)\\
\end{eqnarray*}
Now let us compute ${\mathscr I}(b(\beta))=\sup_{\theta}\left(\theta b(\beta) -\Lambda_\beta(\theta) \right)$. Writing ${A_1+A_2\theta}= \theta{\sigma^2}+(1-\theta){\tau^2}$, the optimizing $\theta$ satisfies
\[b(\beta)=(1-\beta)\left(\log\frac{\sigma}{\tau}-\frac{\frac{1}{2}A_2}{A_1+A_2\theta}\right),\]
which can be solved, giving
\begin{equation*}\theta=-\frac{\frac{1}{2}(1-\beta)}{b(\beta)-(1-\beta)\log(\sigma/\tau)}-\frac{\tau^2}{\sigma^2-\tau^2},\end{equation*}
so that $b(\beta)$ can be evaluated numerically from
\begin{equation}\label{b_variance}
\gamma=(1-\beta)\left(-\frac{1}{2}-\frac{\tau^2}{\sigma^2-\tau^2}\left(\frac{b(\beta)}{1-\beta}-\log\frac{\sigma}{\tau}\right)-\frac{1}{2}\log\left(
\frac{-2\tau^2}{\sigma^2-\tau^2}\left(\frac{b(\beta)}{1-\beta}-\log\frac{\sigma}{\tau}\right)\right)\right).
\end{equation}

\subsection{Test 3: change in scale for dependent data}\label{change_scale}

We now consider the more general situation in which the typical deviations of the process are inflated by a factor $f$. This type of change has applications in the context of communication networks; for details we refer to \cite{KEM2014}. 
More specifically, we concentrate on the case we have that after time $n\beta$ the mean $\bar\mu$ changes into $f\bar\mu$, while the covariance matrix becomes $f^2\Sigma$.
Again, we can shift space so that the first $n\beta$ entries of the alternative mean $\bsn$ equal 0 and the last $n(1-\beta)$ equal $\bar\nu=f\bar\mu-\bar\mu.$
We suppose that ${\bs X}$ corresponds to a stationary sequence of random variables
with possibly `weak dependence' (as defined in \cite[Ch.\ IV]{BROD}); {\sc arma}($p,q$) processes fall in this class. In this section, we assume that the change is introduced abruptly. By this we mean that the memory of observations is not kept after the change which thus results in a new stationary process that is independent from the process before the change.
Because of this, the statistic ${\mathscr L}_{n,\beta}(\bs X)$ of (\ref{LL}) becomes ${\mathscr L}_{n,\beta}(\check{\bs X})=\log\big[ g_{n,\beta}(\check{\bs X})/f_n(\check{\bs X})\big]$, where $\check{\bs X}= (X_{n\beta+1},\ldots, X_n)$. This, using the notation of Section \ref{S3}, reduces to
\begin{eqnarray*}
\lefteqn{\hspace{0cm}\frac{1}{2}\log |\Sigma_{n(1-\beta)}| -\frac{1}{2}\log f^{2n(1-\beta)}|\Sigma_{n(1-\beta)}|
+\frac{1}{2}\check\bsX\TT\Sigma_{n(1-\beta)}^{-1}\check\bsX}\\
&&\hspace{2cm}-\frac{1}{2f^2}(\check\bsX-\bsn_{n(1-\beta)})\TT \Sigma_{n(1-\beta)}^{-1}(\check\bsX-\bsn_{n(1-\beta)})\\
&=&-n(1-\beta)\log f +\frac{1}{2}\check\bsX\TT\Sigma_{n(1-\beta)}^{-1}\check\bsX-\frac{1}{2f^2}(\check\bsX-\bsn_{n(1-\beta)})\TT \Sigma_{n(1-\beta)}^{-1}(\check\bsX-\bsn_{n(1-\beta)}).
\end{eqnarray*}
Using (\ref{momentgeneratingfunction}), it is not hard to verify that the moment generating function ${\mathbb E}_0 \exp (\theta {\mathscr L}_{n,\beta}(\check{\bs X}))$ of our test statistic equals
\[f^{-\theta(1-\beta)n}\left(\sqrt{\theta/f^2+(1-\theta)}\right)^{-(1-\beta)n}\times
\exp\left(-\frac{\theta s_{n,\beta}}{2f^2}\bar\nu^2  +\frac{\theta^2 s_{n,\beta}}{2(\theta f^2+(1-\theta)f^4)}\bar\nu^2\right),\]
with
\[s_{n,\beta}:=\sum_{i=n\beta+1}^n\sum_{i=n\beta+1}^n (\Sigma^{-1})_{i,j},\]
where we recall that  $s_{n,\beta}$ is essentially linear in $n$ and thus the limiting log-moment generating function exists.
The standard machinery now enables us to derive ${b}(\beta)$.

A simplification can be made in case $\bar\nu =0$. This situation occurs when there is no change in mean, while the covariance matrix is multiplied by $f^2$. Then ${b}(\beta)$ follows from
\[\gamma={\mathscr I}(b(\beta))=\sup_\theta\left(\theta {b}(\beta) +\theta(1-\beta)\log f +\frac{1-\beta}{2}\log\left(\frac{\theta}{f^2}+(1-\theta)\right)\right).\]
The optimizing $\theta$ is 
\[-\left(\frac{\frac{1}{2}(1-\beta)}{b(\beta)+(1-\beta)\log f}+\frac{1}{1/f^2-1}\right),\]
so that $b(\beta)$ can be evaluated numerically from
\begin{equation}\label{b_scale}
\gamma=(1-\beta)\left(-\frac{1}{2}-\frac{1}{1/f^2-1}\left(\frac{b(\beta)}{1-\beta}+\log f\right)-\frac{1}{2}\log\left(
\frac{-2}{1/f^2-1}\left(\frac{b(\beta)}{1-\beta}+\log f\right)\right)\right).
\end{equation}
Note that the last equation of Section \ref{change_variance_independent} follows directly from the above equation when $f$ is replaced by $\tau/\sigma$.

\section{Numerical evaluation}\label{S5}

In Section \ref{S4} we have developed changepoint detection tests for dependent sequences. In this section, we evaluate the performance of our proposed method. To this end, we perform a number of simulation experiments. This set-up facilitates evaluating the sensitivity of the procedure, as it enables us to assess its performance in a broad range of scenarios, both in terms of the underlying model, and in terms of the type of change that has taken place in the sequence of observations (in relation to the type of change the sequence is tested against).


We start by explaining the `basic experiment', various variations of which are studied throughout this section. In the basic experiment we simulate an {\sc arma} process with a change from mean 0 to mean 3 and apply the changepoint detection test of Section \ref{SUB41}. (A numerical evaluation for the change in scale test of Section \ref{change_scale} was carried out in \cite{KEM2014}.) 
More specifically, in the basic experiment we carry out the following procedure:

\vb

\begin{itemize}

\item[$\RHD$] In every run we simulate a stationary {\sc ar}(1) or {\sc ma}(1) time series of length 200 that obeys the recursion given in (\ref{RE}) with mean $c=0$ up to observation 99 and mean $c=3$ afterwards, thus having a changepoint at observation 100. The standard deviation $\sigma$ of the $\ve_i$ is set to 1.\footnote{In this experiment ---~consistent with the assumptions in Section \ref{SUB41}~--- the memory $X_{100-1},\ve_{100-1}$ is used as the initial condition for the observation after the change. The transition from the original to the changed process is therefore smooth ---~as opposed to the abrupt change assumed in Section \ref{change_scale}.}

\item[$\RHD$]  We then consider windows of size 50 that we shift along the time series and we test each window for a change in mean. Thus, for window 1 we test observations 1 up to 50 for a changepoint, for window 2 we test observations 2 up to 51 for a changepoint and we continue this procedure up to window 151 which consists of observations 151 up to 200. Note that the first window in which the changepoint is contained is window number 51.

\item[$\RHD$] In order to test for a change in mean within a certain window, we determine whether Inequality (\ref{testcrit}) holds true. To this end, first, the test statistic ${\mathscr L}_{50,\beta}(\bsX)=\log \left[ g_{50,\beta}({\bsX})/ f_{50}({\bs X})\right]$ is computed according to (\ref{LL}). Here $\nu_i$ is 0 for $i<100$ and $\nu_i$ is 3 for $i\geq100$, the covariance matrix $\Sigma=T$ of an {\sc arma} process is computed using the algorithm developed in \cite{MCLE} and $\bsX$ is simulated as described above. Second, the threshold function $b(\beta)$ is computed using (\ref{crit_func_AR1}) for an \textsc{ar}(1) and (\ref{crit_func_MA1}) for an {\sc ma}(1) process.  The significance level $\alpha$ is put to 0.01, so that $\gamma$ in these equations can be found from $e^{-50\,\gamma}=0.01$. Third, we calculate $\frac{1}{50}{\mathscr L}_{50,\beta}(\bsX)-{b}(\beta)$ for  $\beta=\frac{i}{50},i=0,\dots,49$. If the maximum of this difference (taken over $\beta$) is bigger than zero, we raise an alarm. Otherwise we conclude that there is no changepoint in the current window. We repeat this step for all windows. All the steps above are repeated 300 times.

\item[$\RHD$] As soon as we know for each window whether an alarm is raised or not, the performance of the test is evaluated by the following metrics. 
\begin{itemize}
\item[$\bullet$] For every window number the \emph{alarm ratio} is calculated as the number of alarms for that window in 300 runs divided by 300. Note that the alarm ratio for the windows 1 up to 50 gives the \emph{false alarm ratio} per window while for the windows 51 up to 151 it gives the \emph{detection ratio}.
\item[$\bullet$] The \emph{detection delay} is calculated as the time of detection minus the true changepoint. We define the time of detection as the number of the first observation for which we know that a change has happened, that the last observation of the first window in which an alarm was raised after the changepoint occured. For instance, if the changepoint is first detected at time 104 (i.e. the first alarm after the change is raised for window number 55), the delay is 4. We repeat this procedure 300 times, and take the mean of the detection delay over the runs. 
\end{itemize}
\end{itemize}

In the next two sections we discuss the results of the above described experiment, focusing on the alarm ratio in Section \ref{alarm_ratio} and on the detection delay in Section \ref{det_delay}. In Section \ref{exp_variations} we compare the performance of the test for different sizes of the mean shift in order to assess how small of a change in the mean value can be detected. We also examine the sensitivity to the alternative mean chosen in the test setup. We do so by evaluating the performance when testing against a change in mean that is larger than the change we simulate. 

We remark that our straightforward implementation of the procedure in Matlab was executed in 0.1 ms per window. At the same time, in practice a new window will probably be considered only after aggregating a reasonable amount of traffic (which could even be in the order of minutes) in a time bin. In that case 0.1 ms (or even several seconds) of calculation time is fast enough to qualify it as (quasi) on-line. Further improvements can be achieved, for example, by using approximations for the inverse covariance matrix (see, e.g., \cite{Shaman}).

\subsection{Alarm ratio}\label{alarm_ratio}

In this section we analyze the performance of our changepoint detection method by calculating the ratio of (false) alarms as defined above. We will see that for practically relevant coefficients of the {\sc ar}(1) and {\sc ma}(1) processes, the number of false alarms is low. For those coefficients that correspond to a high number of false alarms we explain the reason and describe ways to improve the results.

As examples we consider an {\sc ar} and an {\sc ma} process both with coefficient 0.5, see Figs.\ \ref{DRAR}-- \ref{DRMA}. The dots depict the alarm ratios that we obtained, while the vertical line highlights the earliest window where we could have detected the changepoint. 

\begin{figure}[htbp]
  \centering
  \begin{minipage}[b]{7 cm}
%
%
\begin{tikzpicture}

\begin{axis}[%
width=0.951\figurewidth,
height=\figureheight,
at={(0\figurewidth,0\figureheight)},
scale only axis,
xmin=0,
xmax=150,
xlabel={Window number},
ymin=0,
ymax=1,
ylabel={Alarm ratio},
axis background/.style={fill=white},
axis x line*=bottom,
axis y line*=left
]
\addplot[only marks,mark=o,mark options={},mark size=0.4330pt,color=black] plot table[row sep=crcr]{%
1	0.02\\
2	0.0166666666666667\\
3	0.0166666666666667\\
4	0.0133333333333333\\
5	0.0133333333333333\\
6	0.0233333333333333\\
7	0.0133333333333333\\
8	0.0133333333333333\\
9	0.02\\
10	0.02\\
11	0.0133333333333333\\
12	0.0133333333333333\\
13	0.02\\
14	0.02\\
15	0.0233333333333333\\
16	0.03\\
17	0.02\\
18	0.0133333333333333\\
19	0.02\\
20	0.01\\
21	0.0166666666666667\\
22	0.0133333333333333\\
23	0.0133333333333333\\
24	0.02\\
25	0.0233333333333333\\
26	0.0166666666666667\\
27	0.0333333333333333\\
28	0.02\\
29	0.0266666666666667\\
30	0.0266666666666667\\
31	0.0166666666666667\\
32	0.02\\
33	0.01\\
34	0.00666666666666667\\
35	0.0133333333333333\\
36	0.00666666666666667\\
37	0.02\\
38	0.02\\
39	0.0166666666666667\\
40	0.0366666666666667\\
41	0.02\\
42	0.0166666666666667\\
43	0.0166666666666667\\
44	0.02\\
45	0.01\\
46	0.0133333333333333\\
47	0.0233333333333333\\
48	0.0166666666666667\\
49	0.0133333333333333\\
50	0.0266666666666667\\
51	0.213333333333333\\
52	0.423333333333333\\
53	0.563333333333333\\
54	0.676666666666667\\
55	0.793333333333333\\
56	0.85\\
57	0.906666666666667\\
58	0.943333333333333\\
59	0.96\\
60	0.97\\
61	0.986666666666667\\
62	0.986666666666667\\
63	1\\
64	0.996666666666667\\
65	1\\
66	1\\
67	1\\
68	1\\
69	1\\
70	1\\
71	1\\
72	1\\
73	1\\
74	1\\
75	1\\
76	1\\
77	1\\
78	1\\
79	1\\
80	1\\
81	1\\
82	1\\
83	1\\
84	1\\
85	1\\
86	1\\
87	1\\
88	1\\
89	1\\
90	1\\
91	1\\
92	1\\
93	1\\
94	1\\
95	1\\
96	1\\
97	1\\
98	1\\
99	1\\
100	1\\
101	1\\
102	1\\
103	1\\
104	1\\
105	1\\
106	1\\
107	1\\
108	1\\
109	1\\
110	1\\
111	1\\
112	1\\
113	1\\
114	1\\
115	1\\
116	1\\
117	1\\
118	1\\
119	1\\
120	1\\
121	1\\
122	1\\
123	1\\
124	1\\
125	1\\
126	1\\
127	1\\
128	1\\
129	1\\
130	1\\
131	1\\
132	1\\
133	1\\
134	1\\
135	1\\
136	1\\
137	1\\
138	1\\
139	1\\
140	1\\
141	1\\
142	1\\
143	1\\
144	1\\
145	1\\
146	1\\
147	1\\
148	1\\
149	1\\
150	1\\
};
\addplot [color=black,solid,thick,forget plot]
  table[row sep=crcr]{%
51	0\\
51	1\\
};
\end{axis}
\end{tikzpicture}%
    \caption{\small Alarm ratio per window for an {\sc ar}(1) with coefficient 0.5 and a changepoint at observation 100.}
    \label{DRAR}
  \end{minipage}
  \begin{minipage}[b]{7 cm}
%
%
\begin{tikzpicture}

\begin{axis}[%
width=0.951\figurewidth,
height=\figureheight,
at={(0\figurewidth,0\figureheight)},
scale only axis,
xmin=0,
xmax=150,
xlabel={Window number},
ymin=0,
ymax=1,
ylabel={Alarm ratio},
axis background/.style={fill=white},
title={},
axis x line*=bottom,
axis y line*=left
]
\addplot[only marks,mark=*,mark options={},mark size=0.4330pt,draw=black,fill=black] plot table[row sep=crcr]{%
1	0.0075\\
2	0.005\\
3	0.01\\
4	0.01\\
5	0.0125\\
6	0.0075\\
7	0.0075\\
8	0.0175\\
9	0.015\\
10	0.0225\\
11	0.0125\\
12	0.0075\\
13	0.01\\
14	0.015\\
15	0.0125\\
16	0.015\\
17	0.0175\\
18	0.01\\
19	0.0175\\
20	0.02\\
21	0.02\\
22	0.0175\\
23	0.02\\
24	0.02\\
25	0.0275\\
26	0.025\\
27	0.015\\
28	0.0175\\
29	0.0175\\
30	0.0125\\
31	0.0175\\
32	0.02\\
33	0.01\\
34	0.0125\\
35	0.02\\
36	0.015\\
37	0.02\\
38	0.0125\\
39	0.0075\\
40	0.0175\\
41	0.025\\
42	0.025\\
43	0.0125\\
44	0.025\\
45	0.005\\
46	0.015\\
47	0.0075\\
48	0.015\\
49	0.0125\\
50	0.015\\
51	0.715\\
52	0.7\\
53	0.88\\
54	0.94\\
55	0.9775\\
56	0.985\\
57	0.9975\\
58	1\\
59	1\\
60	1\\
61	1\\
62	1\\
63	1\\
64	1\\
65	1\\
66	1\\
67	1\\
68	1\\
69	1\\
70	1\\
71	1\\
72	1\\
73	1\\
74	1\\
75	1\\
76	1\\
77	1\\
78	1\\
79	1\\
80	1\\
81	1\\
82	1\\
83	1\\
84	1\\
85	1\\
86	1\\
87	1\\
88	1\\
89	1\\
90	1\\
91	1\\
92	1\\
93	1\\
94	1\\
95	1\\
96	1\\
97	1\\
98	1\\
99	1\\
100	1\\
101	1\\
102	1\\
103	1\\
104	1\\
105	1\\
106	1\\
107	1\\
108	1\\
109	1\\
110	1\\
111	1\\
112	1\\
113	1\\
114	1\\
115	1\\
116	1\\
117	1\\
118	1\\
119	1\\
120	1\\
121	1\\
122	1\\
123	1\\
124	1\\
125	1\\
126	1\\
127	1\\
128	1\\
129	1\\
130	1\\
131	1\\
132	1\\
133	1\\
134	1\\
135	1\\
136	1\\
137	1\\
138	1\\
139	1\\
140	1\\
141	1\\
142	1\\
143	1\\
144	1\\
145	1\\
146	1\\
147	1\\
148	1\\
149	1\\
150	1\\
};
\addplot [color=black,solid,thick,forget plot]
  table[row sep=crcr]{%
51	0\\
51	1\\
};
\end{axis}
\end{tikzpicture}%
    \caption{\small Alarm ratio per window for an {\sc ma}(1) with coefficient 0.5 and a changepoint at observation 100.}
    \label{DRMA}
  \end{minipage}
\end{figure}

The picture reveals that we have very few false alarms, their ratio being in the order of 0.01 (as intended since we chose a significance level of 0.01). At the same time, we have achieved the desirable property that the changepoint is detected almost instantly; there is only a small delay. It is noted that {\sc ma}(1) processes fluctuate more frequently than {\sc ar(1)} processes; this may explain the fact that the changepoint is detected earlier for {\sc ma}(1) than for {\sc ar}(1) when both have coefficient 0.5. We come back to the detection delay in Section \ref{det_delay}.


\vb

Above we put the coefficients of the {\sc ma}(1) and {\sc ar}(1) processes equal to 0.5. Now, we want to compare false alarm ratios for a range of different coefficients. To that end we take the mean of the alarm ratios up to the first window where the changepoint is visible; thus, including only windows where every alarm is a false alarm. In this way we obtain Fig.\ \ref{FA}, which shows that for coefficients between $-0.3$ and $0.6$ we obtain an excellent performance in terms of false alarms.
The cases for which the method does not perform well yet can be improved; 
later in this section we point out how the procedure can be adapted to obtain the 
improved curve shown in Fig.\ \ref{FAtuned}. Furthermore, we remark that the proposed method does not have to be used as the only detector but rather can be combined with some other sensors in the effort to reduce the false alarm rate to the acceptable level.

\begin{figure}[htbp]
\centering
\begin{minipage}[b]{7 cm}
%
%
\begin{tikzpicture}

\begin{axis}[%
width=0.951\figurewidth,
height=\figureheight,
at={(0\figurewidth,0\figureheight)},
scale only axis,
xmin=-1,
xmax=1,
xlabel={Coefficient},
ymin=0,
ymax=0.5,
ylabel={False alarm ratio},
axis background/.style={fill=white},
legend style={legend cell align=left,align=left,draw=black},
ytick={0.1,0.2,0.3,0.4,0.5}
]
\addplot [color=black,solid]
  table[row sep=crcr]{%
-0.9	0.125133333333333\\
-0.8	0.0710666666666667\\
-0.7	0.0462\\
-0.6	0.0352666666666667\\
-0.5	0.0265333333333333\\
-0.4	0.0181333333333333\\
-0.3	0.0236666666666667\\
-0.2	0.0144\\
-0.1	0.0164666666666667\\
0	0.0124666666666667\\
0.1	0.0086\\
0.2	0.0115333333333333\\
0.3	0.0096\\
0.4	0.012\\
0.5	0.0172666666666667\\
0.6	0.0322666666666667\\
0.7	0.0582666666666667\\
0.8	0.110666666666667\\
0.9	0.255133333333333\\
};
\addlegendentry{AR(1)};

\addplot [color=black,dashed]
  table[row sep=crcr]{%
-0.9	1\\
-0.8	1\\
-0.7	0.997133333333333\\
-0.6	0.663866666666667\\
-0.5	0.1954\\
-0.4	0.0620666666666667\\
-0.3	0.0363333333333333\\
-0.2	0.0164\\
-0.1	0.0209333333333333\\
0	0.00933333333333333\\
0.1	0.00986666666666667\\
0.2	0.0133333333333333\\
0.3	0.00653333333333333\\
0.4	0.0118666666666667\\
0.5	0.0169333333333333\\
0.6	0.0206666666666667\\
0.7	0.0289333333333333\\
0.8	0.0361333333333333\\
0.9	0.0592666666666667\\
};
\addlegendentry{MA(1)};

\end{axis}
\end{tikzpicture}%
\caption{\small False alarms for a range of different coefficients, basic experiment.}
\label{FA}
\end{minipage}
  \begin{minipage}[b]{7 cm}
%
%
\begin{tikzpicture}

\begin{axis}[%
width=0.951\figurewidth,
height=\figureheight,
at={(0\figurewidth,0\figureheight)},
scale only axis,
xmin=-1,
xmax=1,
xlabel={Coefficient},
ymin=0,
ymax=0.5000000000001,
ylabel={False alarm ratio},
axis background/.style={fill=white},
legend style={legend cell align=left,align=left,draw=black},
ytick={0.1,0.2,0.3,0.4,0.5}
]
\addplot [color=black,solid]
  table[row sep=crcr]{%
-0.9	0.000333333333333333\\
-0.8	0.0004\\
-0.7	6.66666666666667e-05\\
-0.6	0\\
-0.5	6.66666666666667e-05\\
-0.4	0.0002\\
-0.3	0\\
-0.2	0\\
-0.1	0.000266666666666667\\
0	0\\
0.1	0.00106666666666667\\
0.2	0\\
0.3	0\\
0.4	0\\
0.5	0.0002\\
0.6	0.000666666666666667\\
0.7	0.00213333333333333\\
0.8	0.0111333333333333\\
0.9	0.0822666666666667\\
};
\addlegendentry{AR(1)};

\addplot [color=black,dashed]
  table[row sep=crcr]{%
-0.9	1\\
-0.8	0.999933333333333\\
-0.7	0.321333333333333\\
-0.6	0.0217333333333333\\
-0.5	0.0024\\
-0.4	0.0016\\
-0.3	0\\
-0.2	0.0004\\
-0.1	0\\
0	0\\
0.1	0\\
0.2	6.66666666666667e-05\\
0.3	0\\
0.4	6.66666666666667e-05\\
0.5	0.00173333333333333\\
0.6	0\\
0.7	0\\
0.8	0.000133333333333333\\
0.9	0.0006\\
};
\addlegendentry{MA(1)};

\end{axis}
\end{tikzpicture}%
\caption{\small False alarms for a range of different coefficients, adjusted experiment.}
\label{FAtuned}
\end{minipage}
\end{figure}

We now provide an intuitive explanation as to {\it why} our testing procedure tends to perform inadequately for specific parameter values, as we observed in Fig.\ \ref{FA}.
It turns out that the limiting value of $t_{n,\beta}/(n(1-\beta))$, as given  in Lemma 1, is approached slowly for negative coefficients, especially when $\beta$ is big. This effect is illustrated in Figs.\ \ref{LAR}--\ref{LMA} below, where $n$ is plotted against the difference of $t_{n,\beta}/n(1-\beta)$ and the corresponding limit value. As examples we chose a process that showed a good test performance in terms of false alarms (viz.\ an {\sc ar}(1) with coefficient $0.5$) in Fig.\ \ref{LAR}, as well as a process with a very high false alarm rate (viz.\ an {\sc ma}(1) with coefficient $-0.9$) in Fig.\ \ref{LMA}. 

\begin{figure}[htbp]
  \centering
  \begin{minipage}[b]{7 cm}
%
%
\begin{tikzpicture}

\begin{axis}[%
width=0.951\figurewidth,
height=\figureheight,
at={(0\figurewidth,0\figureheight)},
scale only axis,
xmin=10,
xmax=400,
xlabel={n},
ymin=0,
ymax=1,
ylabel={$\text{s}_{\text{n}\beta}\text{/n(1-}\beta\text{) minus  }\tau$},
axis background/.style={fill=white},
legend style={legend cell align=left,align=left,draw=white!15!black}
]
\addplot [color=black,solid]
  table[row sep=crcr]{%
10	0.0833333333333333\\
20	0.0555555555555556\\
30	0.037037037037037\\
40	0.0277777777777778\\
50	0.0222222222222222\\
60	0.0185185185185185\\
70	0.0158730158730159\\
80	0.0138888888888889\\
90	0.0123456790123457\\
100	0.0111111111111111\\
110	0.0101010101010101\\
120	0.00925925925925924\\
130	0.00854700854700857\\
140	0.00793650793650796\\
150	0.00740740740740742\\
160	0.00694444444444442\\
170	0.00653594771241828\\
180	0.00617283950617287\\
190	0.00584795321637427\\
200	0.00555555555555554\\
210	0.00529100529100529\\
220	0.00505050505050503\\
230	0.00483091787439616\\
240	0.00462962962962965\\
250	0.00444444444444442\\
260	0.00427350427350426\\
270	0.00411522633744854\\
280	0.00396825396825395\\
290	0.00383141762452105\\
300	0.00370370370370371\\
310	0.00358422939068098\\
320	0.00347222222222221\\
330	0.00336700336700335\\
340	0.00326797385620914\\
350	0.00317460317460316\\
360	0.00308641975308643\\
370	0.00300300300300299\\
380	0.00292397660818716\\
390	0.00284900284900286\\
400	0.00277777777777777\\
};
\addlegendentry{$\beta\text{=0.1}$};

\addplot [color=black,dashed]
  table[row sep=crcr]{%
1	-0.25\\
2	0.75\\
4	0.5\\
6	0.333333333333333\\
8	0.25\\
10	0.2\\
12	0.166666666666667\\
14	0.142857142857143\\
16	0.125\\
18	0.111111111111111\\
20	0.1\\
22	0.0909090909090909\\
24	0.0833333333333333\\
26	0.0769230769230769\\
28	0.0714285714285715\\
30	0.0666666666666667\\
32	0.0625\\
34	0.0588235294117647\\
36	0.0555555555555556\\
38	0.0526315789473684\\
40	0.05\\
42	0.0476190476190476\\
44	0.0454545454545455\\
46	0.0434782608695652\\
48	0.0416666666666667\\
50	0.04\\
52	0.0384615384615384\\
54	0.037037037037037\\
56	0.0357142857142857\\
58	0.0344827586206897\\
60	0.0333333333333333\\
62	0.0322580645161291\\
64	0.03125\\
66	0.0303030303030303\\
68	0.0294117647058824\\
70	0.0285714285714286\\
72	0.0277777777777778\\
74	0.027027027027027\\
76	0.0263157894736842\\
78	0.0256410256410257\\
80	0.025\\
82	0.024390243902439\\
84	0.0238095238095238\\
86	0.0232558139534884\\
88	0.0227272727272727\\
90	0.0222222222222222\\
92	0.0217391304347826\\
94	0.0212765957446808\\
96	0.0208333333333333\\
98	0.0204081632653061\\
100	0.02\\
102	0.0196078431372549\\
104	0.0192307692307692\\
106	0.0188679245283019\\
108	0.0185185185185185\\
110	0.0181818181818182\\
112	0.0178571428571428\\
114	0.0175438596491228\\
116	0.0172413793103448\\
118	0.0169491525423729\\
120	0.0166666666666667\\
122	0.0163934426229508\\
124	0.0161290322580645\\
126	0.0158730158730159\\
128	0.015625\\
130	0.0153846153846154\\
132	0.0151515151515151\\
134	0.0149253731343283\\
136	0.0147058823529412\\
138	0.0144927536231884\\
140	0.0142857142857143\\
142	0.0140845070422535\\
144	0.0138888888888889\\
146	0.0136986301369863\\
148	0.0135135135135135\\
150	0.0133333333333333\\
152	0.0131578947368421\\
154	0.012987012987013\\
156	0.0128205128205128\\
158	0.0126582278481013\\
160	0.0125\\
162	0.0123456790123457\\
164	0.0121951219512195\\
166	0.0120481927710843\\
168	0.0119047619047619\\
170	0.011764705882353\\
172	0.0116279069767442\\
174	0.0114942528735632\\
176	0.0113636363636364\\
178	0.0112359550561798\\
180	0.0111111111111111\\
182	0.010989010989011\\
184	0.0108695652173913\\
186	0.010752688172043\\
188	0.0106382978723404\\
190	0.0105263157894737\\
192	0.0104166666666667\\
194	0.0103092783505155\\
196	0.0102040816326531\\
198	0.0101010101010101\\
200	0.01\\
202	0.00990099009900991\\
204	0.00980392156862747\\
206	0.00970873786407767\\
208	0.00961538461538464\\
210	0.00952380952380955\\
212	0.00943396226415094\\
214	0.00934579439252337\\
216	0.00925925925925924\\
218	0.00917431192660551\\
220	0.00909090909090909\\
222	0.00900900900900903\\
224	0.00892857142857145\\
226	0.00884955752212391\\
228	0.00877192982456143\\
230	0.00869565217391305\\
232	0.00862068965517243\\
234	0.00854700854700857\\
236	0.00847457627118642\\
238	0.00840336134453784\\
240	0.00833333333333336\\
242	0.00826446280991733\\
244	0.00819672131147542\\
246	0.008130081300813\\
248	0.00806451612903225\\
250	0.00800000000000001\\
252	0.00793650793650796\\
254	0.00787401574803148\\
256	0.0078125\\
258	0.00775193798449614\\
260	0.00769230769230766\\
262	0.00763358778625955\\
264	0.00757575757575757\\
266	0.0075187969924812\\
268	0.0074626865671642\\
270	0.00740740740740742\\
272	0.00735294117647056\\
274	0.00729927007299269\\
276	0.00724637681159418\\
278	0.00719424460431656\\
280	0.00714285714285712\\
282	0.00709219858156029\\
284	0.00704225352112675\\
286	0.00699300699300698\\
288	0.00694444444444442\\
290	0.00689655172413794\\
292	0.00684931506849318\\
294	0.00680272108843538\\
296	0.00675675675675674\\
298	0.00671140939597314\\
300	0.00666666666666665\\
302	0.00662251655629137\\
304	0.00657894736842107\\
306	0.00653594771241828\\
308	0.0064935064935065\\
310	0.00645161290322582\\
312	0.00641025641025639\\
314	0.00636942675159236\\
316	0.00632911392405061\\
318	0.00628930817610063\\
320	0.00624999999999998\\
322	0.00621118012422361\\
324	0.00617283950617287\\
326	0.00613496932515339\\
328	0.00609756097560976\\
330	0.00606060606060604\\
332	0.00602409638554219\\
334	0.0059880239520958\\
336	0.00595238095238093\\
338	0.00591715976331358\\
340	0.00588235294117645\\
342	0.00584795321637427\\
344	0.0058139534883721\\
346	0.00578034682080925\\
348	0.0057471264367816\\
350	0.00571428571428573\\
352	0.00568181818181818\\
354	0.0056497175141243\\
356	0.0056179775280899\\
358	0.00558659217877094\\
360	0.00555555555555554\\
362	0.00552486187845302\\
364	0.00549450549450547\\
366	0.00546448087431695\\
368	0.00543478260869568\\
370	0.00540540540540541\\
372	0.0053763440860215\\
374	0.00534759358288772\\
376	0.00531914893617019\\
378	0.00529100529100529\\
380	0.00526315789473686\\
382	0.00523560209424084\\
384	0.00520833333333331\\
386	0.00518134715025909\\
388	0.00515463917525771\\
390	0.00512820512820511\\
392	0.00510204081632654\\
394	0.00507614213197971\\
396	0.00505050505050503\\
398	0.00502512562814073\\
400	0.005\\
};
\addlegendentry{$\beta\text{=0.5}$};

\addplot [color=black,dotted]
  table[row sep=crcr]{%
1	-0.25\\
10	1\\
20	0.5\\
30	0.333333333333333\\
40	0.25\\
50	0.2\\
60	0.166666666666667\\
70	0.142857142857143\\
80	0.125\\
90	0.111111111111111\\
100	0.1\\
110	0.0909090909090909\\
120	0.0833333333333334\\
130	0.076923076923077\\
140	0.0714285714285715\\
150	0.0666666666666668\\
160	0.0625000000000001\\
170	0.0588235294117648\\
180	0.0555555555555556\\
190	0.0526315789473685\\
200	0.05\\
210	0.0476190476190477\\
220	0.0454545454545455\\
230	0.0434782608695652\\
240	0.0416666666666667\\
250	0.0400000000000001\\
260	0.0384615384615385\\
270	0.0370370370370371\\
280	0.0357142857142858\\
290	0.0344827586206897\\
300	0.0333333333333334\\
310	0.0322580645161291\\
320	0.0312500000000001\\
330	0.0303030303030304\\
340	0.0294117647058824\\
350	0.0285714285714286\\
360	0.0277777777777778\\
370	0.0270270270270271\\
380	0.0263157894736842\\
390	0.0256410256410257\\
400	0.025\\
};
\addlegendentry{$\beta\text{=0.9}$};

\end{axis}
\end{tikzpicture}%
    \caption{\small Difference of $t_{n,\beta}/(n(1-\beta))$ and $\mathscr{T}$ for an {\sc ar}(1) with coefficient $0.5$.}
    \label{LAR}
  \end{minipage}
  \begin{minipage}[b]{7 cm}
%
%
\begin{tikzpicture}

\begin{axis}[%
width=0.951\figurewidth,
height=\figureheight,
at={(0\figurewidth,0\figureheight)},
scale only axis,
xmin=10,
xmax=400,
xlabel={n},
ymin=-100,
ymax=0,
ylabel={$\text{s}_{\text{n}\beta}\text{/n(1-}\beta\text{) minus  }\tau$},
axis background/.style={fill=white},
legend style={at={(0.97,0.03)},anchor=south east,legend cell align=left,align=left,draw=white!15!black}
]
\addplot [color=black,solid]
  table[row sep=crcr]{%
1	-100\\
10	-88.0388823260529\\
20	-68.8919374674708\\
30	-52.1878159301859\\
40	-40.194549371583\\
50	-31.9482914868899\\
60	-26.2064118019626\\
70	-22.0861625467856\\
80	-19.0308376670003\\
90	-16.6951821664562\\
100	-14.8614649926175\\
110	-13.3884971375285\\
120	-12.181910968714\\
130	-11.1767710548435\\
140	-10.3271895314061\\
150	-9.59995717506577\\
160	-8.97052637100214\\
170	-8.4204077134037\\
180	-7.93543949825845\\
190	-7.50461051582883\\
200	-7.11924229469767\\
210	-6.77241021490022\\
220	-6.4585266735465\\
230	-6.17303627298176\\
240	-5.91218977835295\\
250	-5.67287432760401\\
260	-5.45248438564427\\
270	-5.24882259452093\\
280	-5.06002282255868\\
290	-4.88448987947407\\
300	-4.72085187192755\\
310	-4.56792223796015\\
320	-4.4246692588776\\
330	-4.29019139633355\\
340	-4.16369720332033\\
350	-4.04448885338103\\
360	-3.93194855228457\\
370	-3.82552726144633\\
380	-3.7247352872221\\
390	-3.62913438537412\\
400	-3.53833110307151\\
};
\addlegendentry{$\beta\text{=0.1}$};

\addplot [color=black,dashed]
  table[row sep=crcr]{%
1	-100\\
2	-97.8021978021978\\
4	-96.3254196815357\\
6	-94.4707246019884\\
8	-92.3159887804535\\
10	-89.9316810747748\\
12	-87.3844786108474\\
14	-84.7349700736839\\
16	-82.0355733888273\\
18	-79.3295289642457\\
20	-76.6509104061881\\
22	-74.0253598724353\\
24	-71.471230626496\\
26	-69.000876977145\\
28	-66.621913179071\\
30	-64.3383384912424\\
32	-62.1514834024233\\
34	-60.0607702706456\\
36	-58.0643034316973\\
38	-56.1593138352544\\
40	-54.3424857765471\\
42	-52.6101916223185\\
44	-50.9586568327976\\
46	-49.3840734296405\\
48	-47.8826760991226\\
50	-46.4507916909632\\
52	-45.0848700781156\\
54	-43.7815021575094\\
56	-42.5374291152208\\
58	-41.3495458540384\\
60	-40.2149005926173\\
62	-39.130692011648\\
64	-38.0942648772374\\
66	-37.1031047631066\\
68	-36.1548322820074\\
70	-35.2471970939467\\
72	-34.3780718633559\\
74	-33.5454462742883\\
76	-32.74742117153\\
78	-31.9822028688986\\
80	-31.2480976489777\\
82	-30.5435064677357\\
84	-29.8669198706634\\
86	-29.2169131227631\\
88	-28.5921415519784\\
90	-27.9913361038719\\
92	-27.4132991041552\\
94	-26.8569002248242\\
96	-26.3210726490226\\
98	-25.8048094292632\\
100	-25.3071600332451\\
102	-24.8272270711882\\
104	-24.3641631983528\\
106	-23.9171681862201\\
108	-23.4854861556735\\
110	-23.0684029654335\\
112	-22.6652437489667\\
114	-22.2753705931\\
116	-21.8981803516248\\
118	-21.5331025872679\\
120	-21.1795976355319\\
122	-20.8371547840575\\
124	-20.5052905613402\\
126	-20.1835471288279\\
128	-19.8714907706355\\
130	-19.5687104753366\\
132	-19.2748166045181\\
134	-18.9894396430211\\
136	-18.712229026022\\
138	-18.4428520383475\\
140	-18.180992781645\\
142	-17.9263512052619\\
144	-17.6786421969081\\
146	-17.4375947293934\\
148	-17.2029510599407\\
150	-16.9744659787785\\
152	-16.7519061039087\\
154	-16.5350492191298\\
156	-16.3236836525762\\
158	-16.1176076931966\\
160	-15.9166290427581\\
162	-15.720564301112\\
164	-15.5292384825989\\
166	-15.342484561608\\
168	-15.1601430454297\\
170	-14.9820615726612\\
172	-14.8080945355383\\
174	-14.6381027246698\\
176	-14.4719529947504\\
178	-14.3095179499246\\
180	-14.1506756475531\\
182	-13.9953093192236\\
184	-13.8433071079175\\
186	-13.6945618203183\\
188	-13.5489706933128\\
190	-13.4064351737998\\
192	-13.2668607109773\\
194	-13.1301565603353\\
196	-12.9962355986304\\
198	-12.8650141491667\\
200	-12.7364118167507\\
202	-12.6103513317312\\
204	-12.4867584025704\\
206	-12.36556157643\\
208	-12.2466921072912\\
210	-12.130083831155\\
212	-12.0156730479008\\
214	-11.9033984094088\\
216	-11.7932008135745\\
218	-11.68502330387\\
220	-11.578810974127\\
222	-11.4745108782387\\
224	-11.3720719444931\\
226	-11.2714448942755\\
228	-11.1725821648849\\
230	-11.0754378362342\\
232	-10.9799675612114\\
234	-10.8861284994976\\
236	-10.7938792546457\\
238	-10.7031798142418\\
240	-10.6139914929744\\
242	-10.5262768784549\\
244	-10.4399997796371\\
246	-10.3551251776931\\
248	-10.2716191792164\\
250	-10.1894489716227\\
252	-10.1085827806342\\
254	-10.0289898297348\\
256	-9.95064030149334\\
258	-9.8735053006552\\
260	-9.79755681891223\\
262	-9.72276770126166\\
264	-9.64911161387349\\
266	-9.57656301338918\\
268	-9.50509711757658\\
270	-9.43468987727664\\
272	-9.36531794957142\\
274	-9.29695867211785\\
276	-9.2295900385852\\
278	-9.1631906751432\\
280	-9.09773981795036\\
282	-9.03321729159123\\
284	-8.96960348841949\\
286	-8.90687934875979\\
288	-8.84502634193066\\
290	-8.78402644804723\\
292	-8.72386214056709\\
294	-8.66451636954538\\
296	-8.60597254556414\\
298	-8.54821452430724\\
300	-8.49122659174833\\
302	-8.43499344992493\\
304	-8.37950020327263\\
306	-8.32473234549239\\
308	-8.27067574692828\\
310	-8.21731664243167\\
312	-8.16464161969135\\
314	-8.11263760800837\\
316	-8.06129186749511\\
318	-8.01059197868206\\
320	-7.96052583251277\\
322	-7.91108162071113\\
324	-7.8622478265042\\
326	-7.81401321568691\\
328	-7.76636682801168\\
330	-7.71929796889192\\
332	-7.6727962014042\\
334	-7.62685133857745\\
336	-7.58145343595734\\
338	-7.53659278443423\\
340	-7.49225990332407\\
342	-7.44844553369141\\
344	-7.40514063190555\\
346	-7.36233636341986\\
348	-7.32002409676501\\
350	-7.27819539774889\\
352	-7.2368420238528\\
354	-7.19595591881884\\
356	-7.15552920741845\\
358	-7.11555419039705\\
360	-7.0760233395868\\
362	-7.03692929318069\\
364	-6.99826485116361\\
366	-6.96002297089206\\
368	-6.92219676281829\\
370	-6.88477948635332\\
372	-6.84776454586276\\
374	-6.81114548679179\\
376	-6.77491599191339\\
378	-6.73906987769521\\
380	-6.70360109078209\\
382	-6.66850370458845\\
384	-6.6337719159968\\
386	-6.59940004216008\\
388	-6.56538251740133\\
390	-6.53171389020989\\
392	-6.498388820329\\
394	-6.4654020759324\\
396	-6.43274853088622\\
398	-6.40042316209333\\
400	-6.3684210469177\\
};
\addlegendentry{$\beta\text{=0.5}$};

\addplot [color=black,dotted]
  table[row sep=crcr]{%
1	-100\\
10	-95.6223167828945\\
20	-94.0750671490609\\
30	-92.0848424160047\\
40	-89.8647928777945\\
50	-87.4993954354797\\
60	-85.0462833610317\\
70	-82.5486455629712\\
80	-80.03938310629\\
90	-77.5436744326898\\
100	-75.0807476766311\\
110	-72.6651581041489\\
120	-70.3077478095383\\
130	-68.0163909442474\\
140	-65.7965852045286\\
150	-63.6519275213069\\
160	-61.5844995202755\\
170	-59.5951811429\\
180	-57.68390630109\\
190	-55.8498713707153\\
200	-54.0917051175819\\
210	-52.4076069840726\\
220	-50.7954593719864\\
230	-49.2529185328245\\
240	-47.7774878536695\\
250	-46.3665766589892\\
260	-45.0175471032302\\
270	-43.7277512814644\\
280	-42.4945603168066\\
290	-41.3153868791599\\
300	-40.1877023383628\\
310	-39.1090495466352\\
320	-38.0770520727231\\
330	-37.0894205671212\\
340	-36.1439568191225\\
350	-35.2385559680157\\
360	-34.3712072490824\\
370	-33.5399935872864\\
380	-32.7430902953339\\
390	-31.9787630861776\\
400	-31.245365571396\\
};
\addlegendentry{$\beta\text{=0.9}$};

\end{axis}
\end{tikzpicture}%
    \caption{\small  Difference of $t_{n,\beta}/(n(1-\beta))$ and $\mathscr{T}$ for an {\sc ma}(1) with coefficient $-0.9$.}
    \label{LMA}
  \end{minipage}
\end{figure}

We conclude from Figs.\ \ref{LAR}--\ref{LMA} that for the negatively correlated {\sc ma} process we are still far away from the limiting value when $n$ is 400, while for the {\sc ar} process the limiting value is approximated reasonably well already when $n$ is 50 (which corresponds to the chosen window size of 50). 

\vb

In case we do want to handle processes with a high negative correlation we can improve the false alarm rate by adapting our procedure as described in the following paragraphs. As a leading example we consider an {\sc ma}(1) process with coefficient $-0.6$ (see Fig.\ \ref{DRMA05}). One obvious possibility to control the number of false alarms is to lower the significance level $\alpha$ (see Fig.\ \ref{loweralpha}).

\begin{figure}[htbp]
  \centering
  \begin{minipage}[b]{7 cm}
%
%
\begin{tikzpicture}

\begin{axis}[%
width=0.951\figurewidth,
height=\figureheight,
at={(0\figurewidth,0\figureheight)},
scale only axis,
xmin=0,
xmax=150,
xlabel={Window number},
ymin=0,
ymax=1,
ylabel={Alarm ratio},
axis background/.style={fill=white},
title={},
axis x line*=bottom,
axis y line*=left
]
\addplot[only marks,mark=*,mark options={},mark size=0.4330pt,draw=black,fill=black] plot table[row sep=crcr]{%
1	0.6675\\
2	0.6375\\
3	0.645\\
4	0.6375\\
5	0.665\\
6	0.6425\\
7	0.6275\\
8	0.64\\
9	0.67\\
10	0.6675\\
11	0.7075\\
12	0.71\\
13	0.675\\
14	0.6525\\
15	0.66\\
16	0.675\\
17	0.69\\
18	0.6575\\
19	0.7025\\
20	0.6925\\
21	0.69\\
22	0.68\\
23	0.6375\\
24	0.6475\\
25	0.6575\\
26	0.6675\\
27	0.595\\
28	0.6325\\
29	0.71\\
30	0.66\\
31	0.6425\\
32	0.645\\
33	0.6375\\
34	0.65\\
35	0.6375\\
36	0.6775\\
37	0.665\\
38	0.6725\\
39	0.6475\\
40	0.6375\\
41	0.635\\
42	0.6875\\
43	0.69\\
44	0.64\\
45	0.6725\\
46	0.66\\
47	0.675\\
48	0.675\\
49	0.7025\\
50	0.6775\\
51	1\\
52	1\\
53	1\\
54	1\\
55	1\\
56	1\\
57	1\\
58	1\\
59	1\\
60	1\\
61	1\\
62	1\\
63	1\\
64	1\\
65	1\\
66	1\\
67	1\\
68	1\\
69	1\\
70	1\\
71	1\\
72	1\\
73	1\\
74	1\\
75	1\\
76	1\\
77	1\\
78	1\\
79	1\\
80	1\\
81	1\\
82	1\\
83	1\\
84	1\\
85	1\\
86	1\\
87	1\\
88	1\\
89	1\\
90	1\\
91	1\\
92	1\\
93	1\\
94	1\\
95	1\\
96	1\\
97	1\\
98	1\\
99	1\\
100	1\\
101	1\\
102	1\\
103	1\\
104	1\\
105	1\\
106	1\\
107	1\\
108	1\\
109	1\\
110	1\\
111	1\\
112	1\\
113	1\\
114	1\\
115	1\\
116	1\\
117	1\\
118	1\\
119	1\\
120	1\\
121	1\\
122	1\\
123	1\\
124	1\\
125	1\\
126	1\\
127	1\\
128	1\\
129	1\\
130	1\\
131	1\\
132	1\\
133	1\\
134	1\\
135	1\\
136	1\\
137	1\\
138	1\\
139	1\\
140	1\\
141	1\\
142	1\\
143	1\\
144	1\\
145	1\\
146	1\\
147	1\\
148	1\\
149	1\\
150	1\\
};
\addplot [color=black,solid,thick,forget plot]
  table[row sep=crcr]{%
51	0\\
51	1\\
};
\end{axis}
\end{tikzpicture}%
    \caption{\small Alarm ratio per window for an {\sc ma}(1) with coefficient $-0.6$, $\alpha=0.01$.}
    \label{DRMA05}
  \end{minipage}
  \begin{minipage}[b]{7 cm}
%
%
\begin{tikzpicture}

\begin{axis}[%
width=0.951\figurewidth,
height=\figureheight,
at={(0\figurewidth,0\figureheight)},
scale only axis,
xmin=0,
xmax=150,
xlabel={Window number},
ymin=0,
ymax=1,
ylabel={Alarm ratio},
axis background/.style={fill=white},
axis x line*=bottom,
axis y line*=left
]
\addplot[only marks,mark=o,mark options={},mark size=0.4330pt,color=black] plot table[row sep=crcr]{%
1	0.115\\
2	0.115\\
3	0.165\\
4	0.18\\
5	0.125\\
6	0.1\\
7	0.13\\
8	0.12\\
9	0.11\\
10	0.15\\
11	0.14\\
12	0.14\\
13	0.15\\
14	0.105\\
15	0.13\\
16	0.115\\
17	0.14\\
18	0.12\\
19	0.135\\
20	0.125\\
21	0.15\\
22	0.105\\
23	0.09\\
24	0.065\\
25	0.085\\
26	0.1\\
27	0.065\\
28	0.095\\
29	0.11\\
30	0.1\\
31	0.065\\
32	0.12\\
33	0.135\\
34	0.13\\
35	0.105\\
36	0.1\\
37	0.14\\
38	0.13\\
39	0.14\\
40	0.165\\
41	0.16\\
42	0.13\\
43	0.13\\
44	0.135\\
45	0.16\\
46	0.1\\
47	0.13\\
48	0.14\\
49	0.12\\
50	0.115\\
51	0.925\\
52	1\\
53	1\\
54	1\\
55	1\\
56	1\\
57	1\\
58	1\\
59	1\\
60	1\\
61	1\\
62	1\\
63	1\\
64	1\\
65	1\\
66	1\\
67	1\\
68	1\\
69	1\\
70	1\\
71	1\\
72	1\\
73	1\\
74	1\\
75	1\\
76	1\\
77	1\\
78	1\\
79	1\\
80	1\\
81	1\\
82	1\\
83	1\\
84	1\\
85	1\\
86	1\\
87	1\\
88	1\\
89	1\\
90	1\\
91	1\\
92	1\\
93	1\\
94	1\\
95	1\\
96	1\\
97	1\\
98	1\\
99	1\\
100	1\\
101	1\\
102	1\\
103	1\\
104	1\\
105	1\\
106	1\\
107	1\\
108	1\\
109	1\\
110	1\\
111	1\\
112	1\\
113	1\\
114	1\\
115	1\\
116	1\\
117	1\\
118	1\\
119	1\\
120	1\\
121	1\\
122	1\\
123	1\\
124	1\\
125	1\\
126	1\\
127	1\\
128	1\\
129	1\\
130	1\\
131	1\\
132	1\\
133	1\\
134	1\\
135	1\\
136	1\\
137	1\\
138	1\\
139	1\\
140	1\\
141	1\\
142	1\\
143	1\\
144	1\\
145	1\\
146	1\\
147	1\\
148	1\\
149	1\\
150	1\\
151	1\\
};
\addplot [color=black,solid,thick,forget plot]
  table[row sep=crcr]{%
51	0\\
51	1\\
};
\end{axis}
\end{tikzpicture}%
    \caption{\small  Alarm ratio per window for an {\sc ma}(1) with coefficient $-0.6$, $\alpha=0.0001$.}
    \label{loweralpha}
  \end{minipage}
\end{figure}

We can further improve the performance of our testing procedure in terms of false alarms by using a concept similar to the `tuning procedure' proposed in \cite[Section 5]{MAND2}. The main idea behind it is the following. We observed that most false alarms were raised because of a suspected changepoint at the {\it end} of the window, that is, for large $\beta$. (This problem is well known for {\sc llr} tests, see \cite{CZOR}). A simple method to reduce the false alarm rate substantially is to ignore changepoints that correspond to $\beta$ larger than, say, $0.95$ (see Fig.\ \ref{DRtuned}); we  call this adaptation `tuning'. Note that even though we observed that most false alarms occur at the end of the window, tuning also neglects `real' changepoints if they correspond to $\beta>0.95$, and can therefore cause a delayed detection. However, the graph indicates that in the case of an {\sc ma(1)} with coefficient $-0.6$ this approach works remarkably well.

Fig.\ \ref{DRbigwin} shows that we obtain an even better result if we in addition increase the window size to 100.\footnote{To account for the larger window size, in this figure the length of the time series is 300 and the change takes place at time 150.}

\begin{figure}[htbp]
  \begin{minipage}[b]{7 cm}
%
%
\begin{tikzpicture}

\begin{axis}[%
width=0.951\figurewidth,
height=\figureheight,
at={(0\figurewidth,0\figureheight)},
scale only axis,
xmin=0,
xmax=150,
xlabel={Window number},
ymin=0,
ymax=1,
ylabel={Alarm ratio},
axis background/.style={fill=white},
title={},
axis x line*=bottom,
axis y line*=left
]
\addplot[only marks,mark=*,mark options={},mark size=0.4330pt,draw=black,fill=black] plot table[row sep=crcr]{%
1	0.05\\
2	0.06\\
3	0.0575\\
4	0.06\\
5	0.05\\
6	0.0775\\
7	0.0625\\
8	0.0675\\
9	0.06\\
10	0.0575\\
11	0.07\\
12	0.0625\\
13	0.07\\
14	0.0825\\
15	0.0775\\
16	0.0775\\
17	0.0675\\
18	0.0675\\
19	0.0625\\
20	0.055\\
21	0.055\\
22	0.0575\\
23	0.0525\\
24	0.0725\\
25	0.075\\
26	0.0625\\
27	0.06\\
28	0.055\\
29	0.065\\
30	0.0675\\
31	0.06\\
32	0.075\\
33	0.04\\
34	0.055\\
35	0.0575\\
36	0.0625\\
37	0.07\\
38	0.075\\
39	0.0525\\
40	0.07\\
41	0.0675\\
42	0.08\\
43	0.07\\
44	0.0675\\
45	0.07\\
46	0.055\\
47	0.065\\
48	0.065\\
49	0.0675\\
50	0.065\\
51	0.435\\
52	0.895\\
53	1\\
54	1\\
55	1\\
56	1\\
57	1\\
58	1\\
59	1\\
60	1\\
61	1\\
62	1\\
63	1\\
64	1\\
65	1\\
66	1\\
67	1\\
68	1\\
69	1\\
70	1\\
71	1\\
72	1\\
73	1\\
74	1\\
75	1\\
76	1\\
77	1\\
78	1\\
79	1\\
80	1\\
81	1\\
82	1\\
83	1\\
84	1\\
85	1\\
86	1\\
87	1\\
88	1\\
89	1\\
90	1\\
91	1\\
92	1\\
93	1\\
94	1\\
95	1\\
96	1\\
97	1\\
98	1\\
99	1\\
100	1\\
101	1\\
102	1\\
103	1\\
104	1\\
105	1\\
106	1\\
107	1\\
108	1\\
109	1\\
110	1\\
111	1\\
112	1\\
113	1\\
114	1\\
115	1\\
116	1\\
117	1\\
118	1\\
119	1\\
120	1\\
121	1\\
122	1\\
123	1\\
124	1\\
125	1\\
126	1\\
127	1\\
128	1\\
129	1\\
130	1\\
131	1\\
132	1\\
133	1\\
134	1\\
135	1\\
136	1\\
137	1\\
138	1\\
139	1\\
140	1\\
141	1\\
142	1\\
143	1\\
144	1\\
145	1\\
146	1\\
147	1\\
148	1\\
149	1\\
150	1\\
};
\addplot [color=black,solid,thick,forget plot]
  table[row sep=crcr]{%
51	0\\
51	1\\
};
\end{axis}
\end{tikzpicture}%
    \caption{\small Alarm ratio per window for an {\sc ma}(1) with coefficient $-0.6$, $\alpha=0.0001$, when tuning is applied and the window size is 50.}
    \label{DRtuned}
  \end{minipage}
  \begin{minipage}[b]{7 cm}
%
%
\begin{tikzpicture}

\begin{axis}[%
width=0.951\figurewidth,
height=\figureheight,
at={(0\figurewidth,0\figureheight)},
scale only axis,
xmin=0,
xmax=200,
xlabel={Window number},
ymin=0,
ymax=1,
ylabel={Alarm ratio},
axis background/.style={fill=white},
title={},
axis x line*=bottom,
axis y line*=left
]
\addplot[only marks,mark=*,mark options={},mark size=0.4330pt,draw=black,fill=black] plot table[row sep=crcr]{%
1	0.0175\\
2	0.0175\\
3	0.02\\
4	0.0225\\
5	0.02625\\
6	0.0225\\
7	0.03125\\
8	0.02375\\
9	0.03125\\
10	0.0225\\
11	0.0225\\
12	0.01625\\
13	0.02625\\
14	0.0175\\
15	0.02\\
16	0.02375\\
17	0.02125\\
18	0.02875\\
19	0.02125\\
20	0.02375\\
21	0.02625\\
22	0.0275\\
23	0.01875\\
24	0.02375\\
25	0.02375\\
26	0.02125\\
27	0.02125\\
28	0.02\\
29	0.015\\
30	0.02\\
31	0.0225\\
32	0.02375\\
33	0.02\\
34	0.02\\
35	0.025\\
36	0.02375\\
37	0.0175\\
38	0.0175\\
39	0.02125\\
40	0.01875\\
41	0.015\\
42	0.01875\\
43	0.03\\
44	0.03\\
45	0.02375\\
46	0.0225\\
47	0.02375\\
48	0.03125\\
49	0.0325\\
50	0.0375\\
51	0.12\\
52	0.2825\\
53	0.6275\\
54	0.94875\\
55	1\\
56	1\\
57	1\\
58	1\\
59	1\\
60	1\\
61	1\\
62	1\\
63	1\\
64	1\\
65	1\\
66	1\\
67	1\\
68	1\\
69	1\\
70	1\\
71	1\\
72	1\\
73	1\\
74	1\\
75	1\\
76	1\\
77	1\\
78	1\\
79	1\\
80	1\\
81	1\\
82	1\\
83	1\\
84	1\\
85	1\\
86	1\\
87	1\\
88	1\\
89	1\\
90	1\\
91	1\\
92	1\\
93	1\\
94	1\\
95	1\\
96	1\\
97	1\\
98	1\\
99	1\\
100	1\\
101	1\\
102	1\\
103	1\\
104	1\\
105	1\\
106	1\\
107	1\\
108	1\\
109	1\\
110	1\\
111	1\\
112	1\\
113	1\\
114	1\\
115	1\\
116	1\\
117	1\\
118	1\\
119	1\\
120	1\\
121	1\\
122	1\\
123	1\\
124	1\\
125	1\\
126	1\\
127	1\\
128	1\\
129	1\\
130	1\\
131	1\\
132	1\\
133	1\\
134	1\\
135	1\\
136	1\\
137	1\\
138	1\\
139	1\\
140	1\\
141	1\\
142	1\\
143	1\\
144	1\\
145	1\\
146	1\\
147	1\\
148	1\\
149	1\\
150	1\\
151	1\\
152	1\\
153	1\\
154	1\\
155	1\\
156	1\\
157	1\\
158	1\\
159	1\\
160	1\\
161	1\\
162	1\\
163	1\\
164	1\\
165	1\\
166	1\\
167	1\\
168	1\\
169	1\\
170	1\\
171	1\\
172	1\\
173	1\\
174	1\\
175	1\\
176	1\\
177	1\\
178	1\\
179	1\\
180	1\\
181	1\\
182	1\\
183	1\\
184	1\\
185	1\\
186	1\\
187	1\\
188	1\\
189	1\\
190	1\\
191	1\\
192	1\\
193	1\\
194	1\\
195	1\\
196	1\\
197	1\\
198	1\\
199	1\\
200	1\\
};
\addplot [color=black,solid,thick,forget plot]
  table[row sep=crcr]{%
51	0\\
51	1\\
};
\end{axis}
\end{tikzpicture}%
    \caption{\small Alarm ratio per window for an {\sc ma}(1) with coefficient $-0.6$, $\alpha=0.0001$, tuning is applied, the window size is 100.}
    \label{DRbigwin}
  \end{minipage}
\end{figure}

Using these three adjustments ---~that is: (i) a lower significance level of ~$\alpha=0.0001$, (ii)~application of tuning, and (iii)~a larger window of length 100~--- the false alarm performance is substantially better  for most coefficients; compare Fig.\ \ref{FAtuned} with Fig.\ \ref{FA}. However, for {\sc ma}(1) processes with a very high negative correlation (close to $-1$, that is) the window size of 100 is still too small ---~as can be expected from Fig.\ \ref{LMA}. In all other cases the false alarm rate is now close to zero. 

Note that improving the false alarm rate can lead to a lower detection ratio. However, considering the
alarm ratios after the changepoint in Figs.\ \ref{DRMA05}--\ref{DRbigwin}, it is seen that the negative impact of the above adjustments is minor. In some cases a small additional detection delay is introduced, but we always detect the changepoint even when we apply the adjustments. We will see in Section \ref{det_delay} that the negative impact on the delay is smallest for very negative {\sc ma} coefficients, which is exactly the case in which we have the largest number of false alarms (see Fig.\ \ref{FA}), and hence for which the adjustments are most needed. Of course, these results depend also on the magnitude of the new mean after the changepoint. When the mean after the changepoint is large, the adjustment settings can be applied more generally, because the delay decreases (see Section \ref{exp_variations}).

\subsection{Detection delay}\label{det_delay}

After having evaluated how many false alarms are raised before the change, we now wish to assess how fast a changepoint is detected once it occurred. We will see that the delay is low for most {\sc ar} and {\sc ma} coefficients. When using the adjusted settings (to decrease the false alarm ratio), the delay increases, but is still quite low for negative coefficients and very low for \textsc{ma} processes with a very negative coefficient. However, using the adjusted settings for positively correlated processes, increases the detection delay significantly.

In Fig.\ \ref{LOC} we plot the detection delay, which we define as the difference of the detection time and the true changepoint. We do so for a range of different  coefficients of the {\sc ar} and {\sc ma} processes. For comparison, we have included the delays resulting from testing with a single value threshold that was chosen by simulation in such a way that the false alarm rate (approximately) equals the false alarm rate obtained in Fig.~\ref{FA}. Fig.~\ref{LOC} confirms that the changepoint is detected almost immediately for most coefficients. The larger delay for the experiment with simulation-based threshold indicates that a single value threshold can be inferior to a threshold function. 

Fig.~\ref{LOC} also demonstrates that we detect the changepoint earlier for coefficients that correspond to a higher false alarm ratio. 
A notable exception is the case of an {\sc ar}(1) process with a large positive coefficient where both the false alarm ratio (recall Fig. \ref{FA}) and the detection delay are larger. {\sc ar}(1) processes with a high positive correlation tend to behave rather erratically. Therefore, the change is visible later, and moreover, larger jumps have to be tolerated. As an example we may look at a realization of an {\sc ar}(1) process with coefficient 0.9, with a large change from mean 0 to mean 5 at observation 100. The first alarm after the changepoint is raised at window 56, meaning that we locate the changepoint at observation 105. This delay is in line with Fig.\ \ref{ARexample}; actually, by just looking at the process, {it is not clear where to locate the changepoint}.

\vb

{When using the adjusted settings, we detect the changepoint later (compare Fig. \ref{LOC} to Fig.\ \ref{LOCtuned}). When the mean after the change is 3, in the {\sc ar} case the alarm is raised about 4 up to 5 observations late for negative and small positive coefficients. For bigger {\sc ar} coefficients the delay increases sharply. In case of an {\sc ma} process and a change in mean of 3 we are between 4 and 6 observations late for coefficients larger than $-0.3$. For smaller coefficients, the delay is smaller. In short, the adjusted settings have fewest impact on the detection delay for very negative {\sc ma} coefficients while the impact is high for very positive {\sc ar} coefficients. 

We will see in Section \ref{exp_variations} that when the mean after the change is larger, overall the detection delay decreases and thus the negative impact of using the adjusted settings is smaller.  When exactly to apply the adjusted settings depends on the requirements on the false alarm ratio and the detection delay, which differ from application to application. In general, the settings are suited to {\sc ma} processes with a very negative coefficient and to negatively correlated {\sc ar} processes or positively correlated {\sc ma} processes when the change in mean is large (much larger than the standard deviation). When applying the adjusted settings, one should be aware of an increased detection delay for positively correlated {\sc ar} processes.

\begin{figure}[htbp]
  \centering
  \begin{minipage}[t]{7.5 cm}
  \centering
%
%
\begin{tikzpicture}

\begin{axis}[%
width=0.951\figurewidth,
height=\figureheight,
at={(0\figurewidth,0\figureheight)},
scale only axis,
xmin=-1,
xmax=1,
xlabel={Coefficient},
ymin=0,
ymax=40,
ylabel={Delay},
axis background/.style={fill=white},
legend style={at={(0.03,0.97)},anchor=north west,legend cell align=left,align=left,draw=white!15!black}
]
\addplot [color=black,solid]
  table[row sep=crcr]{%
-0.9	0\\
-0.8	0.002\\
-0.7	0.00600000000000001\\
-0.6	0\\
-0.5	0.004\\
-0.4	0.024\\
-0.3	0.042\\
-0.2	0.0840000000000001\\
-0.1	0.198\\
0	0.296\\
0.1	0.46\\
0.2	0.864\\
0.3	1.182\\
0.4	1.728\\
0.5	2.794\\
0.6	4.344\\
0.7	6.256\\
0.8	12.768\\
0.9	24.332\\
};
\addlegendentry{AR: simulated threshold};

\addplot [color=black,dashed]
  table[row sep=crcr]{%
-0.9	0\\
-0.8	0\\
-0.7	0\\
-0.6	0\\
-0.5	0\\
-0.4	0\\
-0.3	0.016\\
-0.2	0.048\\
-0.1	0.146\\
0	0.296\\
0.1	0.46\\
0.2	1.03\\
0.3	0.978\\
0.4	1.728\\
0.5	2.794\\
0.6	4.344\\
0.7	7.14\\
0.8	15.58\\
0.9	35.094\\
};
\addlegendentry{MA: simulated threshold};

\addplot [color=black,dotted]
  table[row sep=crcr]{%
-0.9	0\\
-0.8	0\\
-0.7	0.002\\
-0.6	0.002\\
-0.5	0.032\\
-0.4	0.0640000000000001\\
-0.3	0.128\\
-0.2	0.234\\
-0.1	0.382\\
0	0.582\\
0.1	0.836\\
0.2	1.266\\
0.3	1.46\\
0.4	2.26\\
0.5	2.956\\
0.6	4.086\\
0.7	5.546\\
0.8	7.996\\
0.9	12.922\\
};
\addlegendentry{AR: LD threshold};

\addplot [color=black,dashdotted]
  table[row sep=crcr]{%
-0.9	0\\
-0.8	0\\
-0.7	0\\
-0.6	0\\
-0.5	0.002\\
-0.4	0.0940000000000001\\
-0.3	0.266\\
-0.2	0.406\\
-0.1	0.496\\
0	0.582\\
0.1	0.63\\
0.2	0.782\\
0.3	0.95\\
0.4	1.046\\
0.5	1.096\\
0.6	1.376\\
0.7	1.53\\
0.8	1.67\\
0.9	1.988\\
};
\addlegendentry{MA: LD threshold};

\end{axis}
\end{tikzpicture}%
    \caption{\small Detection delay, basic experiment, change to mean 3; as well as delays obtained with a simulation-based threshold.}
    \label{LOC}
  \end{minipage}
  \hfill
  \begin{minipage}[t]{7 cm}
  \centering
%
%
\begin{tikzpicture}

\begin{axis}[%
width=0.951\figurewidth,
height=\figureheight,
at={(0\figurewidth,0\figureheight)},
scale only axis,
xmin=-1,
xmax=1,
xlabel={Coefficient},
ymin=0,
ymax=40,
ylabel={Delay},
axis background/.style={fill=white},
legend style={at={(0.03,0.97)},anchor=north west,legend cell align=left,align=left,draw=white!15!black}
]
\addplot [color=black,solid]
  table[row sep=crcr]{%
-0.9	3.57\\
-0.8	3.494\\
-0.7	3.598\\
-0.6	3.658\\
-0.5	3.706\\
-0.4	3.718\\
-0.3	3.954\\
-0.2	4.256\\
-0.1	4.41\\
0	4.606\\
0.1	4.81\\
0.2	5.068\\
0.3	5.75\\
0.4	6.608\\
0.5	8.71\\
0.6	12.262\\
0.7	20.714\\
0.8	43.634\\
0.9	43.3453815261044\\
};
\addlegendentry{AR};

\addplot [color=black,dashed]
  table[row sep=crcr]{%
-0.9	0\\
-0.8	0\\
-0.7	0.024\\
-0.6	1.116\\
-0.5	2.31\\
-0.4	3.158\\
-0.3	3.812\\
-0.2	4.184\\
-0.1	4.404\\
0	4.628\\
0.1	4.712\\
0.2	4.892\\
0.3	4.994\\
0.4	5.124\\
0.5	5.298\\
0.6	5.764\\
0.7	6.012\\
0.8	6.446\\
0.9	6.992\\
};
\addlegendentry{MA};

\end{axis}
\end{tikzpicture}%
    \caption{\small  Detection delay, adjusted experiment, change to mean 3.}
    \label{LOCtuned}
  \end{minipage}
\end{figure}

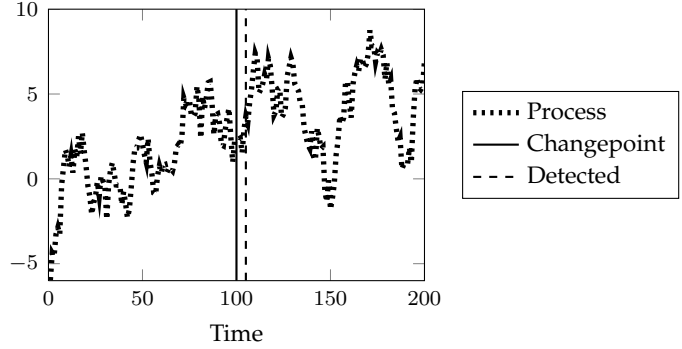
\begin{figure}[h]
\centering
%
%
\begin{tikzpicture}

\begin{axis}[%
width=0.951\figurewidth,
height=\figureheight,
at={(0\figurewidth,0\figureheight)},
scale only axis,
xmin=0,
xmax=200,
xlabel={Time},
ymin=-6,
ymax=10,
axis background/.style={fill=white},
title={},
legend style={at={(1.1,0.5)},anchor=west,legend cell align=left,align=left,draw=black}
]
\addplot [color=black,dotted,line width=2.0pt]
  table[row sep=crcr]{%
1	-5.97138200040203\\
2	-4.37151360198769\\
3	-4.56619758008743\\
4	-3.26692269865116\\
5	-2.8661319840854\\
6	-3.27006903540397\\
7	-0.371388383931727\\
8	0.677313609172501\\
9	1.71888006995439\\
10	1.39810329799805\\
11	0.71673344920511\\
12	1.35182223375812\\
13	0.372647773757855\\
14	1.73594626583203\\
15	1.29352200670572\\
16	2.49773605138861\\
17	1.46976318543525\\
18	2.72225732421195\\
19	1.67877760485192\\
20	0.760803921480223\\
21	-0.437447284488553\\
22	-1.01710844286665\\
23	-2.094639434851\\
24	-1.41094825165793\\
25	-1.05031298508007\\
26	-0.172533827719884\\
27	-0.99497286210185\\
28	0.168905436344737\\
29	-0.126893096655837\\
30	-0.427065444529894\\
31	-2.38088147050246\\
32	-0.329985221106679\\
33	0.522487538128407\\
34	0.956598406874087\\
35	0.283948048422915\\
36	-0.363023396430267\\
37	-0.378345915954797\\
38	-0.3768712510338\\
39	-1.51335154527813\\
40	-0.750140324356533\\
41	-1.44246311592241\\
42	-2.29634084258331\\
43	-1.81098494780126\\
44	-0.545856927623371\\
45	0.910457254322775\\
46	1.81474730659315\\
47	1.81981267609816\\
48	2.29140799402686\\
49	2.63086203038768\\
50	2.02485239061768\\
51	1.43695104496835\\
52	0.974602485835832\\
53	2.37335737787117\\
54	1.00166275358883\\
55	0.315137363777441\\
56	-0.843882105006285\\
57	-0.923100886598835\\
58	0.164027167432862\\
59	0.765153325569737\\
60	0.365425897070301\\
61	0.816654611866443\\
62	0.852058494981386\\
63	0.381595809483139\\
64	0.62089165237418\\
65	0.59541203802253\\
66	0.468537106286814\\
67	0.911012402394989\\
68	1.83598425277548\\
69	2.00532058907248\\
70	1.88689425764091\\
71	2.49590297536027\\
72	4.76247997638537\\
73	4.41216085794458\\
74	4.63077112901652\\
75	4.61585815464988\\
76	3.75527528888926\\
77	3.66026170732907\\
78	3.85439393265211\\
79	5.46491265331754\\
80	4.69489299192277\\
81	2.40596362888482\\
82	4.84754604770549\\
83	4.15126783072913\\
84	4.80628134479811\\
85	5.67264497864528\\
86	5.72859942355307\\
87	4.45709047895486\\
88	2.59625089770426\\
89	3.76671729205105\\
90	3.83046941283756\\
91	3.26781042224048\\
92	1.99968153594117\\
93	1.94312668376687\\
94	2.87793336433774\\
95	3.58631008985533\\
96	4.15803614371367\\
97	2.46754308087103\\
98	0.978216046583191\\
99	2.19654098859676\\
100	2.30297882766584\\
101	2.33910602807638\\
102	2.40425914544263\\
103	1.52965556757887\\
104	2.92011226709407\\
105	3.68425320978873\\
106	3.31297442664344\\
107	4.48672262648788\\
108	5.28185798644526\\
109	6.40448069573574\\
110	7.38959907715858\\
111	7.07655454860562\\
112	5.06499337960028\\
113	5.28662929071744\\
114	5.61804438451624\\
115	6.15387600597592\\
116	5.86214085723664\\
117	7.08556593560026\\
118	6.65642103983785\\
119	5.29092865850876\\
120	3.76120615680979\\
121	3.99594319325438\\
122	5.30557892490421\\
123	5.36552897675481\\
124	3.6214167265175\\
125	3.80163188349792\\
126	4.19150629231492\\
127	6.06493178395891\\
128	6.60509533612792\\
129	7.14369759326151\\
130	6.79262697270081\\
131	5.59801063226578\\
132	5.37151346901473\\
133	5.17510129388338\\
134	4.432232968511\\
135	4.38210758163115\\
136	2.99352586887882\\
137	2.47861268800382\\
138	2.00404424797728\\
139	1.78774602892692\\
140	2.79903270454035\\
141	1.33048054781326\\
142	1.1125506146903\\
143	2.99833523113198\\
144	2.43308466944797\\
145	1.77285023739153\\
146	1.8838487355377\\
147	-0.697432024859033\\
148	-1.34328183855968\\
149	-0.0946449234167988\\
150	-1.79816201325089\\
151	-1.17201039271134\\
152	-0.318692716061417\\
153	1.40583561811352\\
154	2.23138670914722\\
155	3.72105866557729\\
156	3.28408010835253\\
157	3.58384911886579\\
158	5.68716572289269\\
159	5.74422184971869\\
160	4.4589997354391\\
161	3.61631032275281\\
162	5.58860264408019\\
163	6.24887810883329\\
164	6.69286045359217\\
165	6.35760312627519\\
166	7.12740791610293\\
167	7.01388409916137\\
168	6.54761967471531\\
169	6.7429113556794\\
170	7.30688209530148\\
171	8.75774846744177\\
172	7.85073653221078\\
173	7.83698375282848\\
174	6.90369669008313\\
175	7.42816717919765\\
176	7.30386087968764\\
177	7.62514546811535\\
178	7.30295778771863\\
179	5.86154654163584\\
180	6.85191871582599\\
181	5.57250530242625\\
182	5.85881303724551\\
183	4.09514640147616\\
184	3.88425600638999\\
185	2.97212141331909\\
186	2.18694420400732\\
187	2.44646570603505\\
188	0.579030969033795\\
189	0.865797023777958\\
190	0.818154149007177\\
191	0.872724036294387\\
192	1.84352105277403\\
193	0.95116005118684\\
194	3.49759590869228\\
195	5.33152859089584\\
196	5.95378557741031\\
197	5.37453166632042\\
198	5.26134356313781\\
199	5.4408361677088\\
200	6.95814704835256\\
};
\addlegendentry{Process};

\addplot [color=black,thick,solid]
  table[row sep=crcr]{%
100	-6\\
100	10\\
};
\addlegendentry{Changepoint};

\addplot [color=black,thick,dashed]
  table[row sep=crcr]{%
105	-6\\
105	10\\
};
\addlegendentry{Detected};

\end{axis}
\end{tikzpicture}%
\caption{\small Realization of an {\sc ar}(1) process with coefficient 0.9 and a change from mean 0 to 5 at observation 100.}
\label{ARexample}
\end{figure}

\subsection{Sensitivity analysis}\label{exp_variations}
In the above experiments, we chose a shift size $\bar\nu$ and assessed the test's performance for this shift. In the current section we analyze how this performance (in terms of false alarms and detection delay) is affected by the specific value of $\bar\nu$. We will see that -- in accordance with our intuition -- the delay decreases when the change in mean is larger. This may allow us to apply the adjusted settings introduced in Section \ref{alarm_ratio} more generally when the change in mean is large. For the most relevant scenarios (with moderate correlation), the performance in terms of false alarms is good for a broad range of values of $\bar\nu$.

In addition, in our experiments so far, we ran tests in which the mean after the changepoint coincided with the mean we test for. Of course, we would like to have some `robustness'; for that reason we also study in this section the test's performance in case the mean after the changepoint differs from the one that we test for. It turns out that, except for very high positive correlations, the tests are robust against a smaller change than tested for; the detection delay increases slowly when the simulated change becomes smaller.

\vb

\paragraph{\emph{$\RHD$ Varying the size of the change, testing for the mean that we simulated}}
We run the basic experiment, but now we vary the size of the mean shift. Importantly, in these experiments
the mean after the changepoint coincides with the mean we test for. Figs.\ \ref{faarmeans}--\ref{locmameans} describe the tradeoff between an early detection and a low false alarm ratio. As expected, we see that in general it holds that how bigger the change in mean, the smaller the detection delay. The results for the false alarm ratio are somewhat more complicated:

\begin{itemize}

\item For large positive coefficients, we note that the larger the mean the \emph{lower} the number of false alarms. It seems logical that a shift in mean is harder to detect as long as this shift is within the range of the fluctuations typical for the unchanged process. Accordingly, the further $\bar\nu$ exceeds this range the less false alarms we obtain.

\item Surprisingly, for very negative coefficients we see that the opposite: the larger the mean, the \emph{higher} the number of false alarms. For an {\sc ma} process, the false alarm ratio increases much more sharply than for an {\sc ar} process.  To understand this recall that the limit value $\mathscr{T}$ of $t_{n,\beta}/n(1-\beta)$ from Lemma 1 is used to compute the threshold function in (\ref{beetje}). As we saw from Fig. \ref{LMA}, for negative {\sc ma} coefficients $\mathscr{T}$ is substantially \emph{larger} than $t_{n,\beta}/n(1-\beta)$ when $n$ is small. This, in combination with $\bar{\nu}>1$, makes the threshold function more negative than it should be --- the larger $\bar{\nu}$, the more pronounced this effect.

\item When the {\sc ar} or {\sc ma} coefficient is close to zero, neither of the above described effects has a strong impact and the false alarms are systematically low in this case.

\end{itemize}

To summarize, what we have seen is that ---~as we expected~--- detection gets easier as the mean after the change $\bar\nu$ increases. As long as the mean is larger than, say 1 or 1.5 (one or one and halve times the standard deviation of the process), the delay seems acceptable. Concerning the false alarm ratio we have that, for the most relevant case of moderate correlations ({\sc ar} and {\sc ma} coefficients close to zero), the false alarm ratio is low (close to the target of 0.01) for all $\bar\nu$. For highly positively correlated processes the ratio of false alarms is low enough if the change in mean is reasonably large (at least 3, i.e. much larger than the standard deviation of the process). When the correlation is highly negative, the false positive ratio is only low for {\sc ar} processes with a small change in mean (close to the standard deviation). However, the performance of negatively correlated ({\sc ar} with large mean change and {\sc ma}) processes can be improved by using the adjustment settings introduced in Section \ref{alarm_ratio}.

\begin{figure}[h]
  \centering
  \begin{minipage}[b]{7 cm}
%
%
\begin{tikzpicture}

\begin{axis}[%
width=0.951\figurewidth,
height=\figureheight,
at={(0\figurewidth,0\figureheight)},
scale only axis,
xmin=0.5,
xmax=5,
xlabel={Shift size $\bar\nu$},
ymin=0,
ymax=0.300000001,
ylabel={False alarm ratio},
axis background/.style={fill=white},
legend style={legend cell align=left,align=left,draw=black}
]
\addplot [color=black,solid]
  table[row sep=crcr]{%
0.5	0.0102\\
1	0.0174\\
1.5	0.0156\\
2	0.01655\\
2.5	0.02\\
3	0.02715\\
3.5	0.03575\\
4	0.04225\\
4.5	0.05865\\
5	0.08535\\
};
\addlegendentry{$\vr=-0.5$};

\addplot [color=black,dashed]
  table[row sep=crcr]{%
0.5	0.0125\\
1	0.01015\\
1.5	0.0095\\
2	0.01555\\
2.5	0.01405\\
3	0.01655\\
3.5	0.0204\\
4	0.0242\\
4.5	0.019\\
5	0.02535\\
};
\addlegendentry{$\vr=-0.2$};

\addplot [color=black,dotted,line width=2.0pt]
  table[row sep=crcr]{%
0.5	0.02145\\
1	0.0187\\
1.5	0.0157\\
2	0.01865\\
2.5	0.01265\\
3	0.01045\\
3.5	0.00915000000000001\\
4	0.01035\\
4.5	0.0069\\
5	0.00605\\
};
\addlegendentry{$\vr=0.2$};

\addplot [color=black,dashdotted]
  table[row sep=crcr]{%
0.5	0.1443\\
1	0.09445\\
1.5	0.0657\\
2	0.04265\\
2.5	0.03\\
3	0.02025\\
3.5	0.0137\\
4	0.0095\\
4.5	0.00345\\
5	0.0023\\
};
\addlegendentry{$\vr=0.5$};

\end{axis}
\end{tikzpicture}%
    \caption{\small False alarms for different sizes of the mean shift, {\sc ar} case.}
    \label{faarmeans}
  \end{minipage}
  \begin{minipage}[b]{7 cm}
%
%
\begin{tikzpicture}

\begin{axis}[%
width=0.951\figurewidth,
height=\figureheight,
at={(0\figurewidth,0\figureheight)},
scale only axis,
xmin=0.5,
xmax=5,
xlabel={Shift size $\bar\nu$},
ymin=0,
ymax=0.30000001,
ylabel={False alarm ratio},
axis background/.style={fill=white},
legend style={at={(0.03,0.97)},anchor=north west,legend cell align=left,align=left,draw=black}
]
\addplot [color=black,solid]
  table[row sep=crcr]{%
0.5	0.00775\\
1	0.0177\\
1.5	0.03605\\
2	0.06395\\
2.5	0.1106\\
3	0.2119\\
3.5	0.33175\\
4	0.4939\\
4.5	0.662\\
5	0.80035\\
};
\addlegendentry{$\vartheta=-0.5$};

\addplot [color=black,dashed]
  table[row sep=crcr]{%
0.5	0.0128\\
1	0.01085\\
1.5	0.01225\\
2	0.01455\\
2.5	0.0147\\
3	0.0187\\
3.5	0.02055\\
4	0.03005\\
4.5	0.0286\\
5	0.03625\\
};
\addlegendentry{$\vartheta=-0.2$};

\addplot [color=black,dotted,line width=2.0pt]
  table[row sep=crcr]{%
0.5	0.0183\\
1	0.0173\\
1.5	0.01555\\
2	0.01475\\
2.5	0.0106\\
3	0.0126\\
3.5	0.00800000000000001\\
4	0.01135\\
4.5	0.0045\\
5	0.00815\\
};
\addlegendentry{$\vartheta=0.2$};

\addplot [color=black,dashdotted]
  table[row sep=crcr]{%
0.5	0.06685\\
1	0.0477\\
1.5	0.03835\\
2	0.0291\\
2.5	0.0235\\
3	0.01585\\
3.5	0.01325\\
4	0.0065\\
4.5	0.0085\\
5	0.0063\\
};
\addlegendentry{$\vartheta=0.5$};

\end{axis}
\end{tikzpicture}%
    \caption{\small False alarms for different sizes of the mean shift, {\sc ma} case.}
    \label{famameans}
  \end{minipage}
\end{figure}

\begin{figure}[h]
  \centering
  \begin{minipage}[b]{7 cm}
%
%
\begin{tikzpicture}

\begin{axis}[%
width=0.951\figurewidth,
height=\figureheight,
at={(0\figurewidth,0\figureheight)},
scale only axis,
xmin=0.5,
xmax=5,
xlabel={Shift size $\bar\nu$},
ymin=0,
ymax=25,
ylabel={Delay},
axis background/.style={fill=white},
legend style={legend cell align=left,align=left,draw=black}
]
\addplot [color=black,solid]
  table[row sep=crcr]{%
0.5	12.01\\
1	3.10249999999999\\
1.5	1.2525\\
2	0.422499999999999\\
2.5	0.122500000000002\\
3	0.022500000000008\\
3.5	0\\
4	0\\
4.5	0\\
5	0\\
};
\addlegendentry{$\vr=-0.5$};

\addplot [color=black,dashed]
  table[row sep=crcr]{%
0.5	16.415\\
1	4.90000000000001\\
1.5	2.265\\
2	1.0625\\
2.5	0.557500000000005\\
3	0.224999999999994\\
3.5	0.0525000000000091\\
4	0.0250000000000057\\
4.5	0.00499999999999545\\
5	0.00249999999999773\\
};
\addlegendentry{$\vr=-0.2$};

\addplot [color=black,dotted,line width=2.0pt]
  table[row sep=crcr]{%
0.5	21.825\\
1	8.36\\
1.5	4.42\\
2	2.28749999999999\\
2.5	1.55500000000001\\
3	1.1875\\
3.5	0.759999999999991\\
4	0.522500000000008\\
4.5	0.344999999999999\\
5	0.252499999999998\\
};
\addlegendentry{$\vr=0.2$};

\addplot [color=black,dashdotted]
  table[row sep=crcr]{%
0.5	6.95\\
1	5.735\\
1.5	4.91499999999999\\
2	3.52500000000001\\
2.5	3.00999999999999\\
3	2.495\\
3.5	2.00999999999999\\
4	1.78\\
4.5	1.5675\\
5	1.34\\
};
\addlegendentry{$\vr=0.5$};

\end{axis}
\end{tikzpicture}%
    \caption{\small Detection delay for different sizes of the mean shift, {\sc ar} case.}
    \label{locarmeans}
  \end{minipage}
  \begin{minipage}[b]{7 cm}
%
%
\begin{tikzpicture}

\begin{axis}[%
width=0.951\figurewidth,
height=\figureheight,
at={(0\figurewidth,0\figureheight)},
scale only axis,
xmin=0.5,
xmax=5,
xlabel={Shift size $\bar\nu$},
ymin=0,
ymax=25,
ylabel={Delay},
axis background/.style={fill=white},
legend style={legend cell align=left,align=left,draw=black}
]
\addplot [color=black,solid]
  table[row sep=crcr]{%
0.5	8.4725\\
1	2.53999999999999\\
1.5	1.155\\
2	0.557500000000005\\
2.5	0.219999999999999\\
3	0.0250000000000057\\
3.5	0.00249999999999773\\
4	0\\
4.5	0\\
5	0\\
};
\addlegendentry{$\vartheta=-0.5$};

\addplot [color=black,dashed]
  table[row sep=crcr]{%
0.5	16.5825\\
1	5.0975\\
1.5	2.44\\
2	1.27000000000001\\
2.5	0.745000000000005\\
3	0.4375\\
3.5	0.177500000000009\\
4	0.0499999999999972\\
4.5	0.0150000000000006\\
5	0.00499999999999545\\
};
\addlegendentry{$\vartheta=-0.2$};

\addplot [color=black,dotted,line width=2.0pt]
  table[row sep=crcr]{%
0.5	22.2\\
1	7.84999999999999\\
1.5	4.02500000000001\\
2	1.96250000000001\\
2.5	1.1875\\
3	0.677500000000009\\
3.5	0.325000000000003\\
4	0.202500000000001\\
4.5	0.120000000000005\\
5	0.0324999999999989\\
};
\addlegendentry{$\vartheta=0.2$};

\addplot [color=black,dashdotted]
  table[row sep=crcr]{%
0.5	9.9025\\
1	5.67750000000001\\
1.5	3.45750000000001\\
2	1.69\\
2.5	1.22499999999999\\
3	0.537499999999994\\
3.5	0.327500000000001\\
4	0.1875\\
4.5	0.122500000000002\\
5	0.0450000000000017\\
};
\addlegendentry{$\vartheta=0.5$};

\end{axis}
\end{tikzpicture}%
    \caption{\small Detection delay for different sizes of the mean shift, {\sc ma} case.}
    \label{locmameans}
  \end{minipage}
\end{figure} 

\begin{figure}[h]
  \begin{minipage}[t]{7cm}
%
%
\begin{tikzpicture}

\begin{axis}[%
width=0.951\figurewidth,
height=\figureheight,
at={(0\figurewidth,0\figureheight)},
scale only axis,
xmin=0.5,
xmax=5,
xlabel={Shift size (simulated)},
ymin=0,
ymax=50,
ylabel={Delay},
axis background/.style={fill=white},
legend style={legend cell align=left,align=left,draw=black}
]
\addplot [color=black,solid]
  table[row sep=crcr]{%
0.5	3.8775\\
1	0.95999999999998\\
1.5	0.337500000000006\\
2	0.0799999999999841\\
2.5	0.0149999999999864\\
3	0.00499999999999545\\
3.5	0\\
4	0\\
4.5	0\\
5	0\\
};
\addlegendentry{$\vr=-0.5$};

\addplot [color=black,dashed]
  table[row sep=crcr]{%
0.5	11.385\\
1	3.8725\\
1.5	1.50999999999999\\
2	0.657499999999999\\
2.5	0.332499999999982\\
3	0.127499999999998\\
3.5	0.0475000000000136\\
4	0.00749999999999318\\
4.5	0\\
5	0.00249999999999773\\
};
\addlegendentry{$\vr=-0.2$};

\addplot [color=black,dotted]
  table[row sep=crcr]{%
0.5	38.5325\\
1	11.94\\
1.5	5.7175\\
2	3.505\\
2.5	2.25\\
3	1.45749999999998\\
3.5	0.997500000000002\\
4	0.639999999999986\\
4.5	0.402500000000003\\
5	0.215000000000003\\
};
\addlegendentry{$\vr=0.2$};

\addplot [color=black,dashdotted]
  table[row sep=crcr]{%
0.5	88.3125\\
1	34.5525\\
1.5	16.745\\
2	9.17500000000001\\
2.5	6.77500000000001\\
3	4.71250000000001\\
3.5	3.44499999999999\\
4	2.405\\
4.5	1.82999999999998\\
5	1.41500000000002\\
};
\addlegendentry{$\vr=0.5$};

\end{axis}
\end{tikzpicture}%
    \caption{\small Detection delay for different simulated mean shifts in the {\sc ar} case. Always test with shift size $\bar\nu=5$. }
    \label{locarmeanfive}
  \end{minipage}
  \begin{minipage}[t]{7.1 cm}
%
%
\begin{tikzpicture}

\begin{axis}[%
width=0.951\figurewidth,
height=\figureheight,
at={(0\figurewidth,0\figureheight)},
scale only axis,
xmin=0.5,
xmax=5,
xlabel={Shift size (simulated)},
ymin=0,
ymax=50,
ylabel={Delay},
axis background/.style={fill=white},
legend style={legend cell align=left,align=left,draw=black}
]
\addplot [color=black,solid]
  table[row sep=crcr]{%
0.5	0.115000000000009\\
1	0.0349999999999966\\
1.5	0.00749999999999318\\
2	0.00249999999999773\\
2.5	0\\
3	0\\
3.5	0\\
4	0\\
4.5	0\\
5	0\\
};
\addlegendentry{$\vartheta=-0.5$};

\addplot [color=black,dashed]
  table[row sep=crcr]{%
0.5	10.075\\
1	2.8775\\
1.5	1.31\\
2	0.884999999999991\\
2.5	0.467500000000001\\
3	0.257499999999993\\
3.5	0.0949999999999989\\
4	0.0575000000000045\\
4.5	0.0149999999999864\\
5	0.00749999999999318\\
};
\addlegendentry{$\vartheta=-0.2$};

\addplot [color=black,dotted]
  table[row sep=crcr]{%
0.5	33.6225\\
1	10.655\\
1.5	4.72750000000002\\
2	2.48500000000001\\
2.5	1.57749999999999\\
3	0.865000000000009\\
3.5	0.495000000000005\\
4	0.240000000000009\\
4.5	0.0974999999999966\\
5	0.0450000000000159\\
};
\addlegendentry{$\vartheta=0.2$};

\addplot [color=black,dashdotted]
  table[row sep=crcr]{%
0.5	56.97\\
1	17.2675\\
1.5	7.905\\
2	4.0275\\
2.5	2.22750000000002\\
3	1.345\\
3.5	0.680000000000007\\
4	0.332499999999982\\
4.5	0.117500000000007\\
5	0.0250000000000057\\
};
\addlegendentry{$\vartheta=0.5$};

\end{axis}
\end{tikzpicture}%
    \caption{\small Detection delay for different simulated mean shifts in the {\sc ma} case. Always test with shift size $\bar\nu=5$.}
    \label{locmameanfive}
  \end{minipage}
\end{figure}

\paragraph{\emph{$\RHD$ Varying the simulated change in mean, while testing for mean 5}}
We now again vary the simulated mean after the changepoint, but keep the mean that we use in the test setup fixed at 5.

We would expect false alarm rates not to be affected when varying the simulated mean after the changepoint, because false alarms occur before the changepoint. Indeed, we obtain false alarm rates that remain constant for the means we simulated. For coefficients $\geq-0.3$, the false alarm ratio is close to $0.01$, as we aimed for. Consistently with the earlier results, the false alarm ratio is higher for very high coefficients.

We expect the detection delay to increase for a wrongly specified test, where the mean we test for is larger than the actual change. Figs.\ \ref{locarmeanfive}--\ref{locmameanfive} show that the simulated results correspond to this expectation. Nevertheless, it turns out that a change in mean smaller than specified in the test, is tolerated quite well, particularly when the {\sc ar} or {\sc ma} coefficient is small.

\section{Discussion and concluding remarks}

In this paper we have developed \textsc{cusum}-type changepoint detection tests for dependent Gaussian data sequences. The paper includes the setting in which the underlying dataset follows an {\sc arma} structure, a versatile class of models that
has been frequently used to describe traffic streams (and other networking related time series). 
The changepoint tests consist of a log-likelihood test statistic in the spirit of {\sc cusum}, and the corresponding threshold derived from a large-deviations approximation to the false alarm probability. In the literature such \textsc{ld}-based \textsc{cusum}-type tests have so far predominantly focused on procedures for detecting a change in mean in a sequence of independent observations. We have extended the application of this type of test to the case of detecting (1)~a change in mean in correlated normal data, (2)~a change in variance in independent normal data and (3)~a change in scale (that is, the process blows up by a factor) in correlated normal data. Furthermore, the false alarm criterion we employed ensures that the false alarm rate is low for every given window, thus allowing for a low variability of the number of false alarms.

We have demonstrated our changepoint detection test in a number of examples where we tested \textsc{ar}(1) and \textsc{ma}(1) processes against a change in mean. These simulations have shown that the test performs well (in terms of false alarm ratio and detection delay) for \textsc{ar}(1) and \textsc{ma}(1) coefficients between $-0.3$ and $0.6$, as long as the change in mean is larger than the standard deviation of the process. In case of a strong negative correlation or a large change in mean, adaptation of the test settings is possible to further reduce the number of false alarms with minor negative influence on the detection delay. Moreover, the test performance seems to be rather resilient with respect to misspecification of the change size (as used in the test set-up).

\vb

Various next steps could be thought of. A detailed (empirical) comparison to the performance that is achieved under the {\sc arl} criterion is in place. Further, the tests should be modified such that they can be applied to detect a change in the correlation structure within a data sequence. Moreover, other light-tailed distributions may be considered.

\appendix
\section{Proof of Lemma \ref{lemma1}}\label{proof_lemma1}

We first study $v(n):={\mathbb V}{\rm ar} \, S_n$, with $S_n = X_1+\cdots+X_n.$
It follows that
\[S_n - nc =\sum_{i=1}^n \ve_i +\sum_{i=1}^n\sum_{j=1}^p\vr_j\left(X_{i-j}-c\right) +\sum_{i=1}^n\sum_{j=1}^q
\vartheta_j\ve_{i-j}.\]
From this point on we take, without loss of generality, $c=0$.
Recognizing $S_n$ in the right-hand side, bringing all terms involving $S_n$ to the left-hand side, and
taking the variance of both sides, it is now elementary to show that
\begin{equation}
\label{vn}\frac{v(n)}{n} \to \left(\frac{\sigma\left(1+\sum_{j=1}^q \vartheta_j\right)}{1-\sum_{j=1}^p\vr_j}\right)^2;\end{equation}
this identity can alternatively be deduced relying on the spectral density formula for {\sc arma} processes \cite{ROBB}.

Based on `G\"artner-Ellis', with $\pi_n:={\mathbb P}_0\left(S_n \ge n\right),$ 
\[\lim_{n\to\infty}\frac{1}{n}\log \pi_n  = -\frac{1}{2s^2},\] where $s^2$ is the limiting value of $ v(n)/n$ (which we assume to exist).
On the other hand, based on (a discrete-skeleton version of) `Schilder' \cite[Section 4.2]{MAND}, recalling that $T\equiv T_n$ is the covariance matrix of the $X_i$,
\[\lim_{\ve\downarrow 0}\lim_{n\to\infty}\frac{1}{n}\log \pi_n(\ve)=  - \frac{1}{2}\lim_{n\to\infty}
\frac{1}{n}\cdot {\bs 1}T_n^{-1}{\bs 1}=-\frac{1}{2}{\mathscr T}_0,\]
with $\pi_n(\ve):= {\mathbb P}_0\left(\forall i\in\{1,\ldots,n\}:S_i\in(i(1-\ve), i(1+\ve)), S_n\ge n\right).$ We want to prove that

\begin{equation}\label{pies}
\lim_{n\to\infty}\frac{1}{n}\log \pi_n  =\lim_{\ve\downarrow 0}\lim_{n\to\infty}\frac{1}{n}\log \pi_n(\ve),\end{equation}

because if this holds, then the claim of the lemma is an immediate consequence of the fact that $s^{-2}={\mathscr T}_0$. Equation (\ref{pies}) can be proved in three steps.

\begin{itemize}
\item We first observe that, due to `Schilder',
\begin{equation}
\label{pi}\lim_{n\to\infty}\frac{1}{n}\log \pi_n  =\lim_{n\to\infty}\frac{1}{n}\left(-\inf_{{\bs x}\in {\mathscr A}_n} \frac{1}{2}{\bs x}T_n^{-1}{\bs x}\right),\end{equation}
with
${\mathscr A}_n:=\{{\bs x}\,|\,\sum_{i=1}^n x_i \ge n\}.$
It is known \cite[Section 6.1]{MAND} that the optimizing ${\bs x}$, say ${\bs x}^\star$, is such that
\[\sum_{j=1}^i x^\star_j\equiv \sum_{j=1}^i x_j^\star(n) = \frac{{\mathbb C}{\rm ov}\,(S_i,S_n)}{v(n)}\cdot n = \frac{v(n) +v(i)-v(n-i)}{2v(n)}\cdot n.\]
It now follows from (\ref{vn}) that 
\[\lim_{n\to\infty} \sum_{j=1}^i x_j^\star(n) = \lim_{n\to\infty} \frac{ns^2 +is^2-(n-i)s^2}{2ns^2}\cdot n = i.\]
\item
Due to the very same line of reasoning, we also have that
\begin{equation}\label{pie}\lim_{n\to\infty}\frac{1}{n}\log \pi_n(\ve)  =\lim_{n\to\infty}\frac{1}{n}\left(-\inf_{{\bs x}\in {\mathscr B}_n} \frac{1}{2}{\bs x}T_n^{-1}{\bs x}\right),\end{equation}
with, for $\ve>0$,
\[{\mathscr B}_n(\ve):=\left\{{\bs x}\,\left|\, \forall i\in\{1,\ldots,n\}:\sum_{j=1}^ix_j\in(i(1-\ve),i(1+\ve)), 
\sum_{j=1}^n x_j\ge n\right.
\right\}.\]
\item
Obviously, we have that ${\mathscr B}_n(\ve)\subseteq {\mathscr A}_n$ for all $\ve>0$. By construction ${\bs x}^\star$ lies in ${\mathscr A}_n$, but, due to the fact that $\lim_{n\to\infty} \sum_{j=1}^i x_j^\star(n)=i$, we also have that ${\bs x}^\star$ lies in ${\mathscr B}_n(\ve)$
(as $n\to\infty$). As a consequence, Expressions (\ref{pi}) and (\ref{pie}) coincide.
\end{itemize}
Now let $\ve\downarrow 0$, and conclude that
$s^{-2}={\mathscr T}_0$, as claimed.
$\hfill\Box$

\end{document}